\documentclass[journal]{IEEEtran}
\usepackage{cite}
\usepackage{times,amsmath,amssymb,psfrag}
\usepackage{graphicx}
\usepackage{url}
\usepackage{bm}
\usepackage{multicol}
\graphicspath{{./pics/}} \DeclareGraphicsExtensions{.pdf,.jpg,.png}
\usepackage{color}
\usepackage{hyperref}
\usepackage{epsfig}
\usepackage{float}
\usepackage{array}
\usepackage{multirow}
\usepackage{etoolbox}
\usepackage{nomencl}
\usepackage{tikz}
\newcommand{\ct}{\textcolor{black}}

\newcommand\Nom[3][X]{\nomenclature[#1#3]{#2}{#3}}
\renewcommand\nomgroup[1]{%
  \item[\normalsize\itshape\bfseries
  \ifstrequal{#1}{I}{Indices and Sets}{%
  \ifstrequal{#1}{P}{Parameters}{%
  \ifstrequal{#1}{N}{Variables}{%
  \ifstrequal{#1}{X}{Notation for Backward Induction Algorithm}{}}}}]%
}

\makenomenclature
\begin{document}
\title{Dynamic Game-based Maintenance Scheduling of Integrated Electric and Natural Gas Grids with a Bilevel Approach}
\author{
Chong Wang,~\IEEEmembership{Member,~IEEE}, Zhaoyu Wang,~\IEEEmembership{Member,~IEEE}, Yunhe Hou,~\IEEEmembership{Senior Member,~IEEE}, and Kang Ma, ~\IEEEmembership{Member,~IEEE}
\thanks{
This work is supported by the National Science Foundation under CMMI 1745451 and National Natural Science Foundation of China under Grant 51677160.	

C. Wang is with the College of Energy and Electrical Engineering, Hohai University, Nanjing 211100, China (e-mail: wangchonghhu@gmail.com).

Z. Wang is with the Department of Electrical and Computer Engineering, Iowa State University, Ames, IA 5001, USA (email: wzy@iastate.edu).

Y. Hou is with the Department of Electrical \& Electronic Engineering, The University of Hong Kong, Pokfulam, Hong Kong (email: yhhou@eee.hku.hk).

K. Ma is with Electronic \& Electrical Engineering Department, University of Bath, BA2 7AY, UK (e-mail: K.Ma@bath.ac.uk).
 
}
}

\date{}
\maketitle
\begin{abstract}
This paper proposes a dynamic game-based maintenance scheduling mechanism for the asset owners of the natural gas grid and the power grid by using a bilevel approach. In the upper level, the asset owners of the natural gas grid and the power grid schedule maintenance to maximize their own revenues. This level is modeled as a dynamic game problem, which is solved by the backward induction algorithm. In the lower level, the independent system operator (ISO) dispatches the system to minimize the loss of power load and natural gas load in consideration of the system operating conditions under maintenance plans from the asset owners in the upper level. This is modeled as a mixed integer linear programming problem. For the model of the natural gas grid, a piecewise linear approximation associated with the big-M approach is used to transform the original nonlinear model into the mixed integer linear model. Numerical tests on a 6-bus system with a 4-node gas grid and a modified IEEE 118-bus system with a 20-node gas grid show the effectiveness of the proposed model.
\end{abstract}
\IEEEpeerreviewmaketitle
\setlength{\nomitemsep}{0.08cm}
\printnomenclature[2.75cm]
\Nom[I01]{$b,b'$}{Index of power buses.}
\Nom[I02]{$n,n'$}{Index of gas nodes.}
\Nom[I03]{$u$}{Index of generators that will be under maintenance during the given time window.}
\Nom[I04]{$u'$}{Index of generators that will not be under maintenance during the given time window.}
\Nom[I05]{$w$}{Index of gas wells.}
\Nom[I06]{$s$}{Index of gas storages.}
\Nom[I07]{$l$}{Index of lines that will be under maintenance during the given time window.}
\Nom[I08]{$l'$}{Index of lines that will not be under maintenance during the given time window.}
\Nom[I09]{$p$}{Index of pipelines that will be under maintenance during the given time window.}
\Nom[I10]{$p'$}{Index of pipelines that will not be under maintenance during the given time window.}
\Nom[I11]{$pc$}{Index of pipelines with compressors.}
\Nom[I12]{$t,t'$}{Index of time periods.}
\Nom[I12]{$k$}{Index of the linear functions.}
\Nom[I13]{$\Omega _b^U$}{Set of generating units at power bus $b$.}
\Nom[I15]{$\Omega _b^B$}{Set of power buses connected to bus $b$.}
\Nom[I16]{$\Omega _n^W$}{Set of gas wells at natural gas node $n$.}
\Nom[I17]{$\Omega _n^S$}{Set of gas storages at natural gas node $n$.}
\Nom[I18]{$\Omega _n^G$}{Set of gas-fired units at natural gas node $n$.}
\Nom[I20]{$\Omega _n^N$}{Set of gas nodes connected to natural gas node $n$.}

\Nom[P01]{$A_{p,t}^k,A_{p',t}^k$}{Slopes of piecewise linear segments.}
\Nom[P02]{$B_{p,t}^k,B_{p',t}^k$}{Intercepts of piecewise linear segments.}
\Nom[P03]{${B_{b,b'}}$}{Electrical susceptance of line $b$-$b'$.}
\Nom[P04]{$C_{l,t}^M,C_{p,t}^M, C_{u,t}^M$}{Maintenance costs for line $l$,  pipeline $p$ and generator $u$ at $t$.}
\Nom[P05]{$C_{u,t}^F,C_{u,t}^L$}{Fixed cost and linear cost of generator $u$ at $t$.}
\Nom[P05]{$C_{u',t}^F,C_{u',t}^L$}{Fixed cost and linear cost of generator $u'$ at $t$.}
\Nom[P06]{$C_{u',t}^S$}{Start up cost of generator $u$ at $t$.}
\Nom[P07]{$C_{u,t}^P, C_{u',t}^P$}{Revenue of generator $u$ and $u'$ at $t$.}
\Nom[P08]{${C_{w,t}},{C_{s,t}}$}{Cost of gas production and storage at $t$.}
\Nom[P09]{$C_{b,t}^{LS}$}{Cost of power load shedding of bus $b$ at $t$.}
\Nom[P10]{$C_{n,t}^{LS}$}{Cost of gas load shedding of node $n$ at $t$.}
\Nom[P11]{${C_{n,n'}}$}{Weymouth constant of pipeline $n$-$n'$.}
\Nom[P12]{$T_u^{ON},T_u^{OFF}$}{Min on and off time periods of generator $u$.}
\Nom[P13]{$T_w^{ON},T_w^{OFF}$}{Min on and off time periods of gas well $w$.}
\Nom[P14]{$T_l^M,T_p^M,T_u^M$}{Periods of maintenance of line $l$ and pipeline $p$.}
\Nom[P15]{$E{F_u}$}{Efficiency factor of gas-fired unit $u$.}
\Nom[P16]{${\underline{f}_p^k},{\overline{{f}}_p^k}$}{Min/Max gas flow of linear segments for pipeline $p$.}
\Nom[P16]{${\underline{f}_{p'}^k},{\overline{{f}}_{p'}^k}$}{Min/Max gas flow of linear segments for pipeline $p'$.}
\Nom[P17]{${\underline{G}_w}, {\overline{G}_w}$}{Min and max outputs of gas well $w$.}
\Nom[P18]{${\underline{G}_s}, {\overline{G}_s}$}{Min and max outputs of gas storage $s$.}
\Nom[P19]{${L_{n,t}}$}{Gas load of gas node $n$ at $t$ (non-gas-fired unit).}
\Nom[P20]{${L_{b,t}}$}{Power load of bus $b$ at $t$.}
\Nom[P21]{${M^L},{M^P},{M^{\theta}}$}{Large numbers.}

\Nom[P29]{$\overline{N}_t^{L,M}, \overline{N}_t^{P,M}$}{Max number of lines and pipelines under maintenance at $t$.}
\Nom[P29]{$\overline{N}_t^{U,M}$}{Max number of generators under maintenance at $t$.}
\Nom[P30]{$\overline{P}_u, \underline{P}^u$}{Max and min outputs of generator $u$.}
\Nom[P31]{$\overline{P}_{b,b'}^{L}$}{Min/Max capacity of line $b$ - $b'$.}
\Nom[P32]{$\overline{R}_u, \overline{R}_{u'}, \underline{R}_u, \underline{R}_{u'}$}{Ramp-up and ramp-down limits of generator $u$ and $u'$, respectively.}
\Nom[P33]{$\overline{R}_s, \underline{R}_s$}{In-flow and out-flow limits of gas storage $s$.}
\Nom[P34]{$\underline{\theta}_b, \overline{\theta}_b$}{Min/max of phase angles of bus $b$.}
\Nom[P34]{$\overline{\theta}^L$}{Phase angle difference limit.}
\Nom[P35]{$\underline{\pi}_n, \overline{\pi}_n$}{Min/max squared pressures.}
\Nom[P36]{${\lambda _{pc}}$}{Compression factor.}

\Nom[N01]{${F_{n,n',t}}$}{Gas flow from $n$ to $n'$ at $t$.}
\Nom[N02]{$f_{p',t}^k,f_{p,t}^k$}{Piecewise linear gas flow at $t$.}
\Nom[N03]{${G_{w,t}}$}{Gas production of gas well $w$ at $t$.}
\Nom[N04]{${G_{s,t}},{G_{s,t - 1}}$}{Gas inventory of gas storage $s$ at $t$ and $t-1$.}
\Nom[N05]{${L_{u,t}}$}{Gas consumption of gas-fired unit $u$ at $t$.}
\Nom[N06]{$\Delta {L_{n,t}}$}{Gas load shedding of gas node $n$ at $t$.}
\Nom[N07]{$\Delta {L_{b,t}}$}{Power load shedding of bus $b$ at $t$.}
\Nom[N08]{${v_{l,t}},{v_{l,t-1}},{v_{l,t'}}$}{Binary variables to indicate if transmission line $l$ is under maintenance at $t$, $t-1$, and $t'$, respectively. `1' denotes maintenance, otherwise `0'}
\Nom[N09]{${v_{p,t}},{v_{p,t-1}},{v_{p,t'}}$}{Binary variables to indicate if pipeline $p$ is under maintenance at $t$, $t-1$, and $t'$. `1' denotes maintenance, otherwise `0'.}
\Nom[N10]{${v_{u',t}},{v_{u',t - 1}}, {v_{u',t'}}$}{Binary variables to indicate if generator $u'$ is under maintenance at $t$, $t-1$, and $t'$, respectively. `1' denotes maintenance, otherwise `0'.}
\Nom[N10]{${v_{u,t}},{v_{u,t-1}},{v_{u,t'}}$}{Binary variables to indicate the states of generator $u$ at $t$, $t-1$, and $t'$, respectively. `1' and `0' denote on and off-state, respectively. These variables are used for unit commitment.}
\Nom[N11]{${v_{w,t}},{v_{w,t - 1}}, {v_{w,t'}}$}{Binary variables to indicate the states of gas well $w$ at $t$, $t-1$, and $t'$. `1' and `0' denote on and off-state, respectively.}
\Nom[N12]{${P_{u,t}}, {P_{u,t-1}}$}{Power generation of generator $u$ at $t$ and $t-1$, respectively.}
\Nom[N12]{${P_{u',t}}, {P_{u',t-1}}$}{Power generation of generator $u'$ at $t$ and $t-1$, respectively.}
\Nom[N13]{$P_{b,b',t}^L$}{Power from bus $b$ to bus $b'$ at $t$.}
\Nom[N14]{$P{S_{n,t}}, P{S_{n',t}}$}{Pressure of gas nodes $n$ and $n'$ at $t$.}
\Nom[N15]{${o_{u',t}}$}{Binary variable to indicate if generator $u'$ is started up at $t$.}
\Nom[N16]{${\theta _{b,t}},{\theta _{b',t}}$}{Phase angles of buses $b$ and $b'$ at $t$.}
\Nom[N17]{$\eta _{p',t}^k,\eta _{p,t}^k$}{Binary variables to indicate if the piecewise linear function $k$ is selected. `1' denotes `selected', and `0' denotes `non-selected'.}
\Nom[N18]{${\pi _{n,t}},{\pi _{n',t}}$}{Squared pressures of gas nodes $n$ and $n'$.}

\Nom[X01]{$i \in \{ pn,pe\}$}{Player $i$. $pn$ and $pe$ denote the owner of the natural gas grid and the owner of power grid, respectively.}
\Nom[X02]{$e$}{Decision epoch $e$ in the game tree.}
\Nom[X03]{$H_{i,e}$}{Set of decision nodes of the player $i$ at the decision epoch $e$.}
\Nom[X04]{$h_{i,e}$}{A decision node of the player $i$ at the decision epoch $e$, $h_{i,e} \in H_{i,e}$.}
\Nom[X05]{${D_{{h_{i,e}}}}$}{Set of decisions that can be made at the decision  node $h_{i,e}$.}
\Nom[X06]{${d_{{h_{i,e}}}}$}{A decision that can be made at the decision node $h_{i,e}$, ${d_{{h_{i,e}}}} \in {D_{{h_{i,e}}}}$.}
\Nom[X07]{$d_{{h_{i,e}}}^ *$}{Optimal decision at the decision node $h_{i,e}$.}
\Nom[X08]{$\aleph ({d_{{h_{i,e}}}})$}{A decision node followed a decision ${d_{{h_{i,e}}}}$.}
\Nom[X09]{$\nu ({h_{i,e}})$}{A path from the initial decision node to $h_{i,e}$.}
\Nom[X10]{${\zeta ^ * }({h_{i,e}})$}{Optimal path from $h_{i,e}$ to the end decision node.}
\Nom[X11]{${F_i}( \bullet )$}{Function value of player $i$ with a path `$\bullet$'.}
\Nom[X12]{SPNE}{Subgame perfect Nash equilibrium.}

\section{Introduction}
\IEEEPARstart{D}{UE} to higher efficiency, less contamination, and lower costs compared to conventional coal plants, more natural gas-fired units are integrated into power systems. Based on the data from the Energy Information Administration (EIA), 32.1\% of U.S. electricity was supplied by natural gas-fired units in the first month of 2016 \cite{US1}. More than 60\% of new generation required from 2025 to 2040 in the U.S. will be fueled by natural gas \cite{US2}. Growing consumption of natural gas promotes the construction of pipeline networks, which lead to a coupled system associated with the electric grid. Since the natural gas grid is coupled with the electric grid via gas-fired units, it has a great influence on the electric grid \cite{Interdependent1}. An interruption or an outage in the natural gas system may lead to the loss of gas supply for the gas-fired units, which can jeopardize the power system security and result in electric load shedding. For example,  in 2002, a disruption on a single pipeline forced the loss of 2019 MW electricity generation at the Collins generating station near Chicago, and further led to cascading power outages in nearby areas \cite{outage_ex1, outage_ex2}. Therefore, it is necessary to ensure high reliability of the coupled natural gas and electric system. Maintenance scheduling, as an important means to enhance system reliability \cite{Main_Relia1, Resilience_2c}, should include the interactions between the natural gas grid and the electric grid. Particularly, if the assets of the natural gas grid and the electric grid belong to different owners, these owners would like to maximize their own revenues when they make maintenance scheduling. Therefore, it is inevitable to appropriately schedule maintenance for the natural gas grid and the electric grid.

Considering the interdependency between the natural gas grid and the power grid \cite{Chertkov2015541}, the issues such as expansion planning and operation in the coupled grid have been re-investigated \cite{Expansion1, OGF1}. The steady-state models based on Weymouth equations have been applied in studies on optimal unit commitment with security constraints of the natural gas system \cite{UCgas1} and optimal power flow of the electric system with the natural gas system \cite{OPFgas1}. To deal with the nonlinear models of pipelines, a piecewise linear method by using the mixed integer programming formulation can be employed \cite{LinerP1}. In addition, some relaxation-based approaches are proposed to deal with the nonlinear Weymouth equations \cite{Convex1, Convex5, Convex2}. The above studies mainly focus on the steady-state model. To accurately represent the influences of the natural gas grids, the dynamics of the gas pipeline networks are discussed in \cite{DyGas1}.

Currently, there are some studies on maintenance scheduling for the power system. To include the deterioration processes of devices, Markov models are introduced \cite{External1, Chong2_Maintenance, Chong3_Maintenance}. In addition to Markov models, some models based on the mixed integer linear programming (MILP) are proposed to schedule maintenance activities on electric devices \cite{MixMaint1, Chong1_Maintenance}. With the deregulation of the power system, centralized maintenance scheduling is not suitable. Some coordination-based and game-based mechanisms are developed \cite{Deregulation1, game1}. For the natural gas system, condition-based approaches \cite{GasM1}, risk-based approaches \cite{GasM2}, and tree-based approaches \cite{GasM3} are developed to support maintenance scheduling of gas pipelines. However, maintenance scheduling of the integrated natural gas and power systems, especially when the natural gas grid and the power grid are owned by different companies, has been little investigated. {\ct{In practice, the natural gas grid and the power grid may belong to different companies \cite{DifferentC1}. The company of the natural gas grid owns assets in the natural gas grid, and the company of the power grid owns assets in the power grid \cite{EPRG1}. When scheduling maintenance, the owners of the natural gas grid and the electric grid both expect to maximize their revenues \cite{MaxR1}. In addition, there are independent system operators (ISOs) for the coupled system, and they coordinate, control and monitor the operation of the system \cite{FERC1}. For different owners of the natural gas grid and the power grid, they have the rights to schedule maintenance for their own assets, and ISOs are in charge of the system operation to ensure high reliability of the systems. Since the two grids may belong to different owners, different sequences of determining maintenance scheduling for different owners may happen, and the sequence of determining maintenance scheduling can impact the revenues of the owners for the natural gas grid and the electric grid. Therefore, we propose a dynamic game-based maintenance scheduling mechanism for the electric and natural gas grids by using a bilevel approach. In the upper level, the different owners of the natural gas grid and the power grid schedule maintenance with the objective to maximize their own revenues. The interactions between different owners are formulated as a dynamic game problem, which is solved by the backward induction algorithm. In the lower level, ISOs dispatch the natural gas grid and the power grid to minimize the loss of power load and gas load with the operating conditions under the maintenance schedules from the different owners.}} This is formulated as a mixed integer linear programming problem. A piecewise linear approximation {\ct{with the big-M approach}} is used to transform the original nonlinear model of the natural gas grid into a mixed integer linear model. A 6-bus system with a 4-node gas grid and a modified IEEE 118-bus system with a 20-node gas grid are used to verify the effectiveness of the proposed model. 

The main contributions are listed as follows. 1) A dynamic game-based maintenance scheduling mechanism for different owners of the electric grid and the natural gas grid is proposed by using a bilevel model. The dynamic game model addresses the sequence of scheduling maintenance for different owners of the natural gas grid and the power grid. The bilevel model deals with that the asset owners determine the maintenance windows of the assets and ISOs determine the system operating conditions. 2) The subgame perfect Nash equilibrium (SPNE) for the owners of the natural gas grid and the electric grid is obtained by using the backward induction algorithm. 3) The influences of maintenance durations and the number of piecewise linear functions on the SPNE are analyzed. Suggestions on maintenance durations are provided.

The remainder of this paper is organized as follows. Section II shows the maintenance scheduling formulation, including the framework of the proposed model, the upper level model and the lower level model. Section III presents the solution, and the case studies are demonstrated in Section IV. The work is concluded in Section V. 

\section{Maintenance Scheduling Formulation}
This section shows the formulation of the dynamic game-based maintenance scheduling with a bilevel approach. First, the framework of the proposed model is presented. Second, the upper level model, including the model for the owners of the natural gas grid and the power grid, is introduced. Third, the lower level model for the ISO with the operating constraints is established.
\subsection{Framework of the proposed model}
\begin{figure}[!h]
	\centering
	\includegraphics[width=8.5cm]{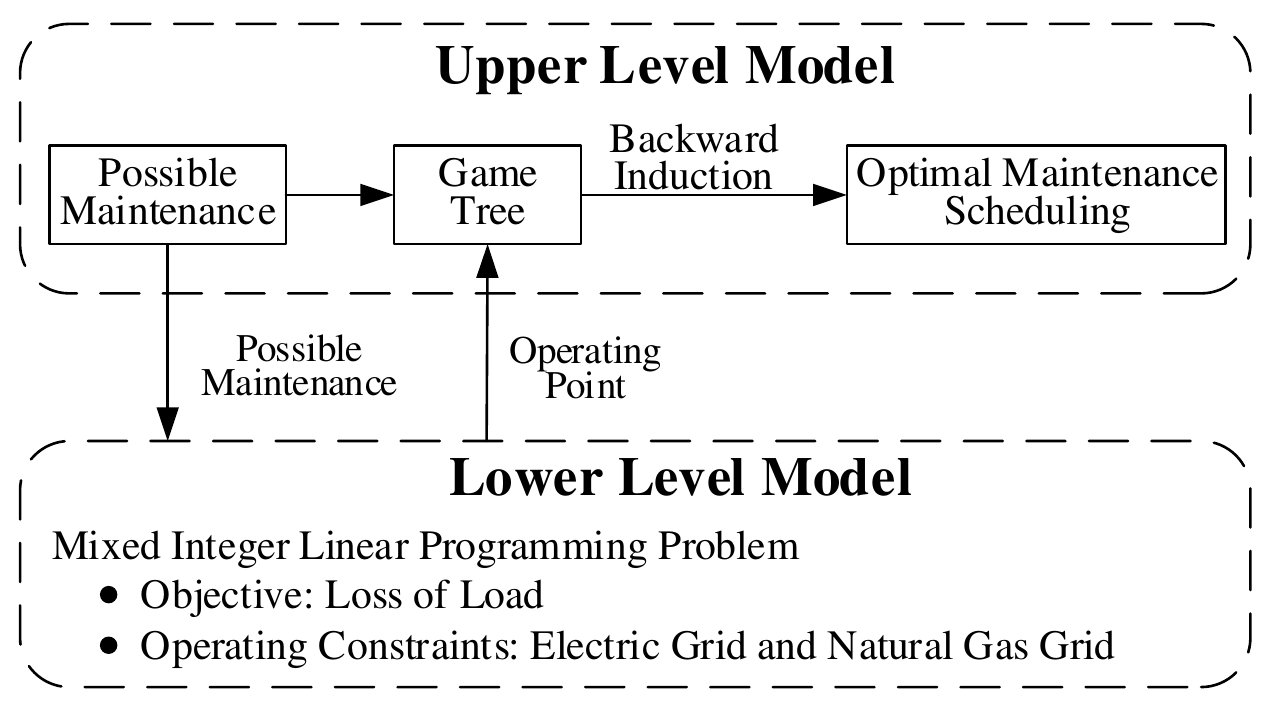}	
	\caption{Framework of the proposed model}	
	\label{framework_fig}
\end{figure}

Fig. \ref{framework_fig} shows the framework of the proposed dynamic game-based model, which is a bilevel model. The upper-level model schedules maintenance with the objective to maximize the revenues for the owners of the natural gas grid and the power grid, and the lower-level model calculates the system operating points for the ISO with the objective to minimize the loss of power load and gas load. In this paper, there are the following assumptions: 1) For the system, there is one owner of the natural gas grid and one owner of the power grid. In addition, there is one ISO who is responsible for dispatching the system. 2) Knowledge about one owner is available to the other owner and the ISO. 

In the upper level, possible maintenance can be listed with the scheduling constraints and the scheduling sequence. For possible maintenance, the lower level model is employed to minimize the loss of load with the operating constraints of the natural gas grid and the power grid, and it is a mixed integer linear programming problem. Based on the operating points, the game tree with payoffs can be established. Then, the game tree is solved by using the backward induction algorithm, and SPNE is obtained. SPNE corresponds to the maintenance scheduling for the owners of the natural gas grid and the power grid.    

\subsection{Upper Level Model}
\subsubsection{Power Grid Owner's Model}
In the upper-level model, the owner of the power grid schedules maintenance activities while keeping a maximum revenue. The objective can be expressed as follows. 
\begin{eqnarray}
\label{Eq.UpperEleObj}
\ \begin{array}{l}
\max \quad \mathop {\underline {\sum\limits_u {\sum\limits_t {\left( {C_{u,t}^P \cdot {P_{u,t}}} \right)} } } }\limits_{1.a}  + \mathop {\underline {\sum\limits_{u'} {\sum\limits_t {\left( {C_{u',t}^P \cdot {P_{u',t}}} \right)} } } }\limits_{1.b}  - \\[20pt]
\;\mathop {\underline {\sum\limits_u {\sum\limits_t {\left( {C_{u,t}^M \cdot {v_{u,t}}} \right)} } } }\limits_{1.c}  - \mathop {\underline {\sum\limits_l {\sum\limits_t {\left( {C_{l,t}^M \cdot {v_{l,t}}} \right)} } } }\limits_{1.d}  - \\[20pt]
\;\mathop {\underline {\sum\limits_u {\sum\limits_t {\left( {C_{u,t}^F \cdot {v_{u,t}} + C_{u,t}^L \cdot {P_{u,t}} } \right)} } } }\limits_{1.e}  - \\[20pt]
\;\mathop {\underline {\sum\limits_{u'} {\sum\limits_t {\left( {C_{u',t}^F \cdot {v_{u',t}} + C_{u',t}^L \cdot {P_{u',t}} + C_{u',t}^S \cdot {o_{u',t}}} \right)} } } }\limits_{1.f}  - \\[20pt]
\;\mathop {\underline {\sum\limits_b {\sum\limits_t {\left( {C_{b,t}^{LS} \cdot \Delta {L_{b,t}}} \right)} } } }\limits_{1.g} 
\end{array}
\end{eqnarray}
s.t.
\begin{align}
\label{Eq.UpperEleCon1}
& \begin{array}{l}
{v_{l,t - 1}} - {v_{l,t}} + {v_{l,t'}} \le 1\\
\quad \quad \quad 1 \le t' - (t - 1) \le T_l^M,\forall l,t
\end{array} \\[1pt]
\label{Eq.UpperEleCon1_1}
& \begin{array}{l} {v_{u,t - 1}} - {v_{u,t}} + {v_{u,t'}} \le 1\quad \quad \\ \quad \quad\quad 1 \le t' - (t - 1) \le T_u^M,\forall u,t \end{array} \\
\label{Eq.UpperEleCon2}
& \sum\limits_t {\left( {1 - {v_{l,t}}} \right) = T_l^M} \quad \quad \forall l \quad \quad\quad \quad \\
\label{Eq.UpperEleCon2_1}
& \sum\limits_t {\left( {1 - {v_{u,t}}} \right) = T_u^M} \quad \quad \forall l \quad \quad\quad \quad \\
\label{Eq.UpperEleCon3}
& \sum\limits_{l} {\left( {1 - {v_{l,t}}} \right) \le \overline{N}_t^{L,M}} \quad \quad \forall t \\
\label{Eq.UpperEleCon3_1}
& \sum\limits_{u} {\left( {1 - {v_{u,t}}} \right) \le \overline{N}_t^{U,M}} \quad \quad \forall t \\
\label{Eq.UpperEleCon4}
& {v_{l,t}}, {v_{u,t}} \in \{ 0,1\} \quad \forall l,t \\
\label{Eq.UpperEleCon5}
&\begin{array}{l}
\Delta {L_{b,t}},{v_{u',t}},{P_{u,t}},\\ \quad \quad {P_{u',t}}, {o_{u',t}} \in \arg \left\{ {{\bf{Lower}}\;{\bf{Level}}\;{\bf{Model}}} \right\} \end{array} 
\end{align}
where (\ref{Eq.UpperEleObj}.a) and (\ref{Eq.UpperEleObj}.b) show the revenues from power generation, (\ref{Eq.UpperEleObj}.c) is the cost of implementing maintenance on generators, and (\ref{Eq.UpperEleObj}.d) is the cost of implementing maintenance on lines. (\ref{Eq.UpperEleObj}.e) is the operating cost of generators that will be under maintenance during the given time window. It is assumed that the generators to be under maintenance will not participate in unit commitment in the lower-level model. (\ref{Eq.UpperEleObj}.f) is the operating cost of generators that will not be under maintenance during the given time window, and these generators can participate in unit commitment in the lower-level model. (\ref{Eq.UpperEleObj}.g) is the penalty cost due to unserved power load.
(\ref{Eq.UpperEleCon1}) is the maintenance duration constraint for the line $l$, and (\ref{Eq.UpperEleCon1_1}) is the maintenance duration constraint for the generator $u$.
(\ref{Eq.UpperEleCon2}) and (\ref{Eq.UpperEleCon2_1}) ensure that maintenance will be implemented during the given time window, (\ref{Eq.UpperEleCon3}) is the constraint of the maximum allowable number of lines under maintenance in one period, and (\ref{Eq.UpperEleCon3_1}) is the constraint of the maximum allowable number of generators under maintenance in one period. (\ref{Eq.UpperEleCon4}) is the binary constraint. (\ref{Eq.UpperEleCon5}) shows the variables optimized in the lower-level model.

In (\ref{Eq.UpperEleObj})-(\ref{Eq.UpperEleCon5}), the variables $v_{l,t}, v_{u,t}$ $\forall l,u,t$, which represent the maintenance states, are determined in the upper-level model. This is because the owner of the power grid only determines maintenance windows for the devices. The variables $\Delta {L_{b,t}},{v_{u',t}},{P_{u,t}}, {P_{u',t}}, {o_{u',t}}$ $\forall u,u',b,t$, which represent the system operating conditions, are determined by the ISO in the lower-level model.

\subsubsection{Natural Gas Grid Owner's Model}
In the upper-level model, the objective of the natural gas grid can be expressed as follows. 
\begin{eqnarray}
\label{Eq.UpperGasObj}
\ \begin{array}{l}
\max \quad \mathop {\underline {\sum\limits_w {\sum\limits_t {\left( {{C_{w,t}} \cdot {G_{w,t}}} \right)} } } }\limits_{10.a}  + \mathop {\underline {\sum\limits_s {\sum\limits_t {\left( {{C_{s,t}} \cdot {G_{s,t}}} \right)} } } }\limits_{10.b}  - \\[20pt]
\quad \quad \;\;\mathop {\underline {\sum\limits_p {\sum\limits_t {\left( {C_{p,t}^M \cdot {v_{p,t}}} \right)} } } }\limits_{10.c}  - \mathop {\underline {\sum\limits_n {\sum\limits_t {\left( {C_{n,t}^{LS} \cdot \Delta {L_{n,t}}} \right)} } } }\limits_{10.d} 
\end{array} 
\end{eqnarray}
s.t.
\begin{align}
\label{Eq.UpperGasCon1}
&\begin{array}{l}
{v_{p,t - 1}} - {v_{p,t}} + {v_{p,t'}} \le 1\\
\quad \quad \quad 1 \le t' - (t - 1) \le T_p^M,\forall p,t
\end{array} \\
\label{Eq.UpperGasCon2}
&\sum\limits_{t} {\left( {1 - {v_{p,t}}} \right) = T_p^M} \;  \forall p \\
\label{Eq.UpperGasCon3}
&\sum\limits_{p} {\left( {1 - {v_{p,t}}} \right) \le \bar N_t^{P,M}} \;  \forall t \\
\label{Eq.UpperGasCon4}
&{v_{p,t}} \in \{ 0,1\} \; \forall p,t \\
\label{Eq.UpperGasCon5}
&\Delta {L_{n,t}},{v_{w,t}},{G_{w,t}},{G_{s,t}} \in \arg \left\{ {{\bf{Lower}}\;{\bf{Level}}\;{\bf{Model}}} \right\} 
\end{align}
where (\ref{Eq.UpperGasObj}.a) and (\ref{Eq.UpperGasObj}.b) are the revenues for the owner of the natural gas grid. (\ref{Eq.UpperGasObj}.c) is the maintenance cost, (\ref{Eq.UpperGasObj}.d) is the penalty cost caused by unserved natural gas load, (\ref{Eq.UpperGasCon1}) is the constraint of  maintenance on pipelines, (\ref{Eq.UpperGasCon2}) guarantees that maintenance will be implemented during the given time window, (\ref{Eq.UpperGasCon3}) enforces the maximum number of pipelines that can be maintained in one period, (\ref{Eq.UpperGasCon4}) is the binary constraint, (\ref{Eq.UpperGasCon5}) shows the variables optimized in the lower-level model.

In (\ref{Eq.UpperGasObj})-(\ref{Eq.UpperGasCon5}), the variables $v_{p,t}$ $\forall p,t$, representing the maintenance states, are determined by the owner of the natural gas grid in the upper-level model. The variables $\Delta {L_{n,t}},{v_{w,t}},{G_{w,t}},{G_{s,t}}$ $\forall n,w,s,t$, representing the operating conditions, are determined by the ISO in the lower-level model.

\subsection{Lower Level Model}
\subsubsection{Objective function}
In the lower-level model, the ISO minimizes the loss of power load and natural gas load while satisfying the operating constraints with the scheduled maintenance from the owners of the natural gas grid and the power grid. The objective of the lower-level model is expressed as follows. 
\begin{eqnarray}
\label{Eq.LowerObj}
\ \min \quad \sum\limits_{b} {\sum\limits_{t} {\left( {C_{b,t}^{LS} \cdot \Delta {L_{b,t}}} \right)} }  + \sum\limits_{n} {\sum\limits_{t} {\left( {C_{n,t}^{LS} \cdot \Delta {L_{n,t}}} \right)} }
\end{eqnarray}
\subsubsection{Constraints}
When minimizing the objective function, the operating constraints of the natural gas grid and the electric grid should be satisfied.  
\paragraph{Natural gas supply constraints}
In each period, the natural gas supply constraints should be satisfied. 
\begin{align}
\label{Eq.LowerCon1}
&\;\;{\underline{G}_w} \cdot {v_{w,t}} \le {G_{w,t}} \le {\overline{G}_w} \cdot {v_{w,t}}\quad \quad \forall w,t \\
\label{Eq.LowerCon2}
&\begin{array}{l}
- {v_{w,t - 1}} + {v_{w,t}} - {v_{w,t'}} \le 0\quad \\
\quad \quad \quad \quad\quad 1 \le t' - (t - 1) \le T_w^{ON},\forall w,t
\end{array} \\
\label{Eq.LowerCon3}
&\begin{array}{l}
{v_{w,t - 1}} - {v_{w,t}} + {v_{w,t'}} \le 1\quad \quad \\
\quad \quad \quad \quad 1 \le t' - (t - 1) \le T_w^{OFF},\forall w,t
\end{array} \\
\label{Eq.LowerCon4}
&\;{\underline{G}_s} \le {G_{s,t}} \le {\overline{G}_s}\quad \forall s,t  \\
\label{Eq.LowerCon5}
&\;{\underline{R}_s} \le {G_{s,t}} - {G_{s,t + 1}} \le {\overline{R}_s}\quad \forall s,t  
\end{align}
where (\ref{Eq.LowerCon1}) shows the limit of gas production of each gas well, (\ref{Eq.LowerCon2}) and (\ref{Eq.LowerCon3}) show the constraints of on-off states of the gas wells, (\ref{Eq.LowerCon4}) enforces the storage levels of the gas storages, (\ref{Eq.LowerCon5}) represents the in-flow and out-flow rates of the gas storages.
\paragraph{Natural gas balance at each gas node}
Natural gas balance at each gas node should be satisfied.
\begin{align}
\label{Eq.LowerCon6}
& \begin{array}{l}
\sum\limits_{w \in \Omega _n^W} {{G_{w,t}}}  + \sum\limits_{s \in \Omega _n^S} {\left( {{G_{s,t - 1}} - {G_{s,t}}} \right)}  + \sum\limits_{n' \in \Omega _n^N} {{F_{n,n',t}}} \\
\quad \quad  - \sum\limits_{u \in \Omega _n^G} {\left( {{L_{u,t}} \cdot E{F_u}} \right)}  - \left( {{L_{n,t}} - \Delta {L_{n,t}}} \right) = 0\quad \forall n,t
\end{array} \\
\label{Eq.LowerCon7}
&\quad 0 \le \Delta {L_{n,t}} \le {L_{n,t}}\quad \quad\quad\quad \quad\quad\quad\quad\; \forall t,n 
\end{align}
where (\ref{Eq.LowerCon6}) enforces the natural gas balance at each node, and (\ref{Eq.LowerCon7}) denotes the constraint for natural gas load shedding.

\begin{figure}[]
	\centering
	\includegraphics[width=7cm]{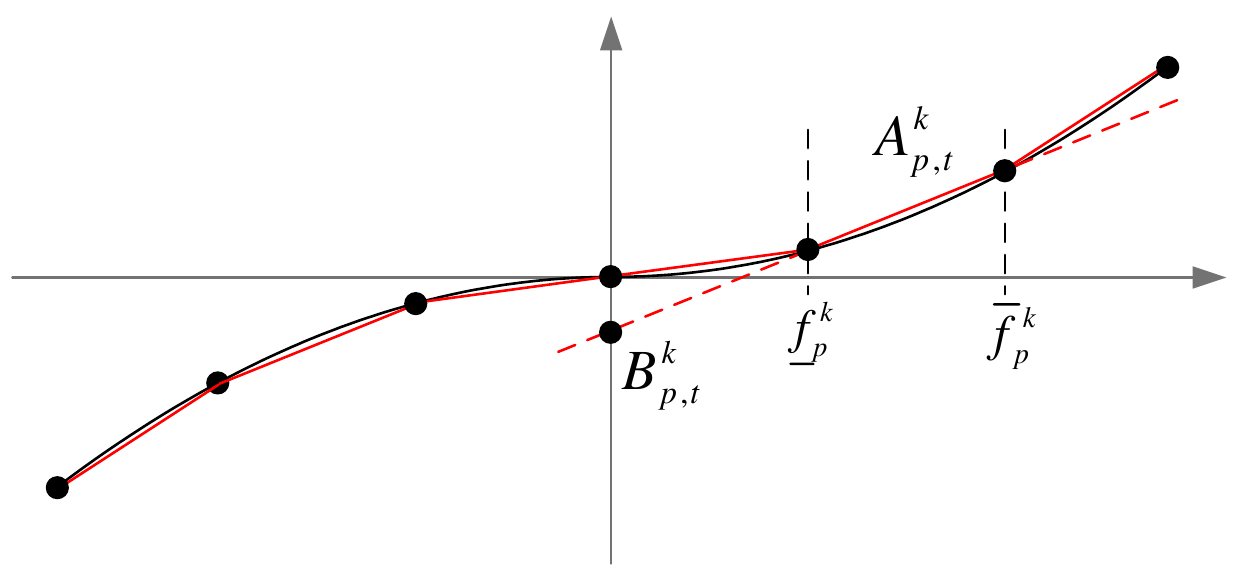}
	\caption{Piecewise linear approximation of pipeline flow.}
	\label{piecewise1_fig}
\end{figure}

\paragraph{Natural gas flow constraints}
{\ct{The gas flow in a pipeline in service can be expressed as the nonlinear Weymouth equation (\ref{Eq.LowerCon8_0}).
\begin{align}
\label{Eq.LowerCon8_0}
& {F_{n,n',t}} \cdot \left| {{F_{n,n',t}}} \right| = C_{n,n'}^2 \cdot (PS_{n,t}^2 - PS_{n',t}^2)\; \forall n,n',t
\end{align}
A piecewise linear formulation by using the mixed integer programming, as shown in Fig. \ref{piecewise1_fig}, is employed to deal with the Weymouth equation.}} Substitute $PS^2$ with $\pi$, and the left side of (\ref{Eq.LowerCon8_0}) can be approximately expressed as the sum of a group of piecewise linear functions. Each linear function is represented by a slope $A$, an axis intercept and a binary variable $\eta$.

For the pipelines that will not be under maintenance during the given time window, the model can be expressed as follows. 
\begin{align}
\label{Eq.LowerCon8}
& \begin{array}{l}
\sum\limits_k {\left( {A_{p',t}^k \cdot f_{p',t}^k + B_{p',t}^k \cdot \eta _{p',t}^k} \right)} \\
\quad \quad \quad  = C_{n,n'}^2\left( {{\pi _{n,t}} - {\pi _{n',t}}} \right)\quad \forall t,p',\;(n,n') \in p'
\end{array} \\
\label{Eq.LowerCon11}
& \sum\limits_k {\eta _{p',t}^k}  = 1\quad \forall t,p' \\
\label{Eq.LowerCon9}
&\;\eta _{p',t}^k \cdot {\underline{f}}_{p'}^k \le f_{p',t}^k \le \eta _{p',t}^k \cdot {\overline{f}}_{p'}^k\quad \forall k,t,p' \\
\label{Eq.LowerCon10}
&\;{F_{n,n',t}} = \sum\limits_k {f_{p',t}^k} \quad \forall t,p',\;(n,n') \in p'  
\end{align}
where the sum of a group of piecewise linear functions in the left side of (\ref{Eq.LowerCon8}) approximates the term ${F_{n,n',t}} \cdot \left| {{F_{n,n',t}}} \right| $ in (\ref{Eq.LowerCon8_0}).
In Fig. \ref{piecewise1_fig}, there are six piecewise linear functions, i.e., $k = \{ 1,2, \cdots ,6\}$. Since the variable ${\eta _{p',t}^k} \in \{0,1\}$, (\ref{Eq.LowerCon11}) ensures that only one linear function will be selected. (\ref{Eq.LowerCon9}) limits the variable $f_{p',t}^k$ for the $k$th linear function. With the constraints (\ref{Eq.LowerCon8}), (\ref{Eq.LowerCon11}) and (\ref{Eq.LowerCon9}), the gas flow in the pipeline $n-n'$ can be expressed as (\ref{Eq.LowerCon10}). 

{\ct{For a pipeline that will be under maintenance during the given time window,  its state (i.e., the variable $v_{p,t}$) is determined in the upper-level model. The piecewise linear approximation model is shown as follows.}}
\begin{align}
\label{Eq.LowerCon12}
& \begin{array}{l}
\sum\limits_k {\left( {A_{p,t}^k \cdot f_{p,t}^k + B_{p,t}^k \cdot \eta _{p,t}^k} \right)}  \le C_{n,n'}^2\left( {{\pi _{n,t}} - {\pi _{n',t}}} \right)\\
\quad \quad \quad \quad \quad  + \left( {1 - {v_{p,t}}} \right) \cdot {M^P}\quad \quad \forall t,p,\;(n,n') \in p
\end{array} \\
\label{Eq.LowerCon13}
& \begin{array}{l}
\sum\limits_k {\left( {A_{p,t}^k \cdot f_{p,t}^k + B_{p,t}^k \cdot \eta _{p,t}^k} \right)}  \ge C_{n,n'}^2\left( {{\pi _{n,t}} - {\pi _{n',t}}} \right)\\
\quad \quad \quad \quad \quad  - \left( {1 - {v_{p,t}}} \right) \cdot {M^P}\quad \quad \forall t,p,\;(n,n') \in p
\end{array} \\
\label{Eq.LowerCon14}
&\;\;\eta _{p,t}^k \cdot \underline f_p^k \le f_{p,t}^k \le \eta _{p,t}^k \cdot \overline f_p^k\quad \forall k,t,p \\
\label{Eq.LowerCon15}
&\;\;{F_{n,n',t}} = \sum\limits_k {f_{p,t}^k} \quad \forall t,p,\;(n,n') \in p \\
\label{Eq.LowerCon16}
&\;\;\sum\limits_k {\eta _{p,t}^k}  = {v_{p,t}}\quad \forall t,p   
\end{align}
where (\ref{Eq.LowerCon12}) and (\ref{Eq.LowerCon13}) represent the constraints of gas flow in pipelines. By introducing the binary variable $v_{p,t}$ and a sufficiently large $M^P$, (\ref{Eq.LowerCon12}) and (\ref{Eq.LowerCon13}) are redundant when the corresponding pipelines are under maintenance at $t$. (\ref{Eq.LowerCon14}) and (\ref{Eq.LowerCon15}) have the similar meanings of (\ref{Eq.LowerCon9}) and (\ref{Eq.LowerCon10}), respectively. (\ref{Eq.LowerCon16}) relates the variable ${\eta _{p,t}^k}$ and ${v_{p,t}}$. For example, when the pipeline $p$ at $t$ is off due to maintenance, i.e., ${v_{p,t}}=0$, ${\eta _{p,t}^k}$ for all piecewise linear functions will be zeros. In this case, (\ref{Eq.LowerCon12})-(\ref{Eq.LowerCon15}) are all satisfied. When the pipeline $p$ at $t$ is on, i.e., ${v_{p,t}}=1$, the model is the same as (\ref{Eq.LowerCon8}) - (\ref{Eq.LowerCon10}). The binary variable $v_{p,t}$ is determined in the upper-level model, and the binary variable ${\eta _{p,t}^k}$ is determined in the lower-level model.

\paragraph{Natural gas pressure constraints}
The pressure of each natural gas node should be within the limit.
\begin{align}
\label{Eq.LowerCon17}
& {\underline{\pi}_n} \le {\pi _{n,t}} \le {\overline{\pi} _n}\quad\quad\quad\quad\quad \forall n,t \\
\label{Eq.LowerCon18}
& {\pi _{n',t}} \le {\lambda _{pc}} \cdot {\pi _{n,t}}\quad\quad\quad\quad \forall pc,t,(n,n') \in pc  
\end{align}
where (\ref{Eq.LowerCon17}) is the constraint for the pressure. (\ref{Eq.LowerCon18}) shows the constraint of pressures between the in-coming gas node and the out-coming gas node of a pipeline with a compressor.

\paragraph{Power flow constraints}
The limits for power flow through lines should be satisfied.
\begin{align}
\label{Eq.LowerCon19}
&\;{B_{b,b'}} \cdot \left( {{\theta _{b,t}} - {\theta _{b',t}}} \right) = P_{b,b',t}^L\quad \forall t,l',\;(b,b') \in l' \\
\label{Eq.LowerCon20}
&\left| {P_{b,b',t}^L} \right| \le \overline P_{b,b'}^L\quad \forall t,l',\;(b,b') \in l' \\
\label{Eq.LowerCon21}
&\begin{array}{l}
{B_{b,b'}} \cdot \left( {{\theta _{b,t}} - {\theta _{b',t}}} \right) - P_{b,b',t}^L + (1 - {v_{l,t}}) \cdot {M^L} \ge 0\\
\quad \quad \quad \quad \quad \quad \quad \quad \quad \quad \quad \quad \quad\quad \quad \forall t,l,\;(b,b') \in l
\end{array} \\
\label{Eq.LowerCon22}
& \begin{array}{l}
{B_{b,b'}} \cdot \left( {{\theta _{b,t}} - {\theta _{b',t}}} \right) - P_{b,b',t}^L - (1 - {v_{l,t}}) \cdot {M^L} \le 0\\
\quad \quad \quad \quad \quad \quad \quad \quad \quad \quad \quad \quad \quad \quad \quad \forall t,l,\;(b,b') \in l
\end{array} \\
\label{Eq.LowerCon23}
&\;\left| {P_{b,b',t}^L} \right| \le \overline P_{b,b'}^L \cdot {v_{l,t}}\quad \quad \forall t,l,\;(b,b') \in l 
\end{align}
where (\ref{Eq.LowerCon19}) shows the relation between voltage angles and power flow through lines that will not be under maintenance during the given time window, (\ref{Eq.LowerCon20}) shows the corresponding limit of power flow. (\ref{Eq.LowerCon21}) and (\ref{Eq.LowerCon22}) represent the physical relations between voltage angles and power flow through lines that will be under maintenance during the given time window. (\ref{Eq.LowerCon23}) is the corresponding limit of power flow.

\paragraph{Generator's state constraint} When scheduling maintenance, it is necessary to include the on-off constraints of generators. 
\begin{align}
\label{Eq.LowerCon24_1}
& \begin{array}{l}
- {v_{u',t - 1}} + {v_{u',t}} - {v_{u',t'}} \le 0\quad \quad \\
\quad \quad \quad \quad\quad \quad 1 \le t' - (t - 1) \le T_u^{ON},\forall u',t
\end{array} \\
\label{Eq.LowerCon25_1}
& \begin{array}{l}
{v_{u',t - 1}} - {v_{u',t}} + {v_{u',t'}} \le 1\\
\quad \quad \quad\quad \quad 1 \le t' - (t - 1) \le T_u^{OFF},\forall u',t
\end{array} \\
\label{Eq.LowerCon26_1}
& - {v_{u',t - 1}} + {v_{u',t}} - {o_{u',t}} \le 0\quad \quad \forall u',t 
\end{align}
where (\ref{Eq.LowerCon24_1}) denotes the minimum on-time constraint of the generator $u'$, (\ref{Eq.LowerCon25_1}) is the minimum off-time constraint of the generator $u'$. These constraints are used for unit commitment. (\ref{Eq.LowerCon26_1}) is the start-up constraint of the generator $u'$. 

\paragraph{Power balance} At each time, the power system should satisfy power balance. 
\begin{align}
\label{Eq.LowerCon27}
&\begin{array}{l}
\sum\limits_{u \in \Omega _b^U} {{P_{u,t}}}  + \sum\limits_{u' \in \Omega _b^U} {{P_{u',t}}}  - \left( {{L_{b,t}} - \Delta {L_{b,t}}} \right) + \\
\quad \quad \quad \quad \quad \quad\quad \quad \quad\quad  \sum\limits_{b' \in \Omega _b^B} {P_{b,b',t}^L}  = 0\quad \quad \forall t,b
\end{array} \\[5pt]
\label{Eq.LowerCon30}
&\;\;\; 0 \le \Delta {L_{b,t}} \le {L_{b,t}}\quad\quad\quad \quad \forall t,b \\
\label{Eq.LowerCon28}
&\;\;\; {P_{u,t}} = {L_{u,t}} \cdot E{F_u}\quad \;\forall u \in \Omega _n^G,n,t 
\end{align}
where (\ref{Eq.LowerCon27}) enforces power balance at each bus in each time period. (\ref{Eq.LowerCon30}) is the constraint of load shedding. (\ref{Eq.LowerCon28}) denotes the relation between natural gas and real power produced from gas-fired units.

\paragraph{Phase angle constraint} The constraint for the phase angle at each bus, and the constraint for the phase angle difference between two buses on each line should be satisfied. 
\begin{align}
\label{Eq.LowerCon29_1}
&\left| {{\theta _{b,t}} - {\theta _{b',t}}} \right| \le {v_{l,t}} \cdot {{\overline \theta^L }} + (1 - {v_{l,t}}) \cdot {M^\theta }\quad \forall t,l,(b,b') \in l \\
\label{Eq.LowerCon29_2}
&\left| {{\theta _{b,t}} - {\theta _{b',t}}} \right| \le {{\overline \theta }^L}\quad \forall t,l',(b,b') \in l'
\end{align}
where the constraint (\ref{Eq.LowerCon29_1}) shows the phase angle difference limit on the line that will be maintained during the given time window, and (\ref{Eq.LowerCon29_2}) shows the phase angle difference limit on the line that will not be maintained during the given time window.

\paragraph{Ramp-up and ramp-down constraints} The generators should satisfy the ramp-up and ramp-down constraints.
\begin{align}
\label{Eq.LowerCon31}
&\;\; {{\underline{P}}_u} \cdot {v_{u,t}} \le {P_{u,t}} \le {{\overline{P}}_u} \cdot {v_{u,t}}\quad \quad \forall u,t \\
\label{Eq.LowerCon32}
& \begin{array}{l}
{P_{u,t}} - {P_{u,t - 1}} \le (2 - {v_{u,t - 1}} - {v_{u,t}}) \cdot {{\underline{P}}_u} + \\
\quad \quad \quad \quad \quad \quad  (1 + {v_{u,t - 1}} - {v_{u,t}}) \cdot {{\overline{R}}_u}\quad \quad \forall u,t
\end{array} \\
\label{Eq.LowerCon33}
& \begin{array}{l}
{P_{u,t - 1}} - {P_{u,t}} \le (2 - {v_{u,t - 1}} - {v_{u,t}}) \cdot {{\underline{P}}_u} + \\
\quad \quad \quad \quad \quad \quad (1 - {v_{u,t - 1}} + {v_{u,t}}) \cdot {{\underline{R}}_u}\quad \quad \forall u,t
\end{array} 
\end{align}
where (\ref{Eq.LowerCon31}) is the capacity limits of the generator $u$, (\ref{Eq.LowerCon32}) and (\ref{Eq.LowerCon33}) are ramp-up and ramp-down constraints of the generators $u$. When $v_{u,t}=1$ and $v_{u,t-1}=0$, the constraints become ${P_{u,t}} - {P_{u,t - 1}} \le  {{\overline{R}}_{u}}$ and ${P_{u,t - 1}} - {P_{u,t}} \le  {{\underline{R}}_{u}}$, and these two constraints are the conventional ramping limits. When $v_{u,t}=1$ and $v_{u,t-1}=0$, the constraints can be rewritten as ${P_{u,t}}  \le  {{\underline{P}}_{u}}$ and $-{P_{u,t}}  \le  {{\underline{P}}_{u}} +2{{\underline{R}}_{u}}$. Associated with (\ref{Eq.LowerCon31}), we can conclude ${P_{u,t}}={{\underline{P}}_{u}}$. This means that the generation of the generator $u$ is ${{\underline{P}}_{u}}$ after it is restarted. When $v_{u,t}=0$ and $v_{u,t-1}=1$, we have ${P_{u,t-1}}={{\underline{P}}_{u}}$ and ${P_{u,t}}=0$. It means the generation of the generator $u$ will be dispatched to ${{\underline{P}}_{u}}$ at $t-1$ when the generator $u$ should be offline at $t$. For the generators that will not be under maintenance during the given time window, the ramp-up and ramp-down constraints are similar with (\ref{Eq.LowerCon31})-(\ref{Eq.LowerCon33}).

\section{Solutions}
This section presents the solution procedures to obtain maintenance plans for the owners of the natural gas grid and the power grid. The solution procedures include 1) establishing the game tree, 2) calculating payoffs of each path in the game tree, and 3) obtaining SPNE by using the backward induction algorithm.

\subsection{Game Tree Establishment}
In the upper-level model, the owners of the natural gas grid and the power grid schedule their own maintenance. With the constraints (\ref{Eq.UpperEleCon1})-(\ref{Eq.UpperEleCon4}) and the constraints (\ref{Eq.UpperGasCon1})-(\ref{Eq.UpperGasCon4}), a game tree representing possible sequential decisions can be established. Since the decision variables in (\ref{Eq.UpperEleCon1})-(\ref{Eq.UpperEleCon4}) and (\ref{Eq.UpperGasCon1})-(\ref{Eq.UpperGasCon4}) are discrete, the established game tree is a finite game tree. 

\begin{figure}[!h]
	\centering
	\includegraphics[width=3.0cm]{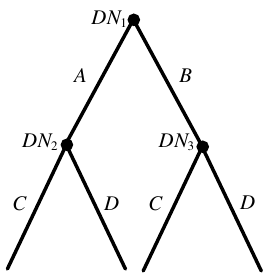}
	\caption{A simplified game tree.}
	\label{gametree_fig}
\end{figure}

For example, the owner of the natural gas grid can make maintenance scheduling $A$ and $B$, and the owner of the power grid can make maintenance scheduling $C$ and $D$. When the owner of the natural gas grid first makes maintenance scheduling, the game tree is illustrated in Fig. \ref{gametree_fig}. $DN_1$ is the decision node for the owner of the natural gas grid. $DN_2$ and $DN_3$ are the decision nodes for the owner of the power grid.

\subsection{Payoff Calculation}
For each path in the game tree, i.e., the scheduled maintenance, the ISO dispatches the system to ensure the high reliability and minimizes the loss of power load and natural gas load. The variables $v_{l,t}, v_{u,t}, v_{p,t}$ in each path come from the upper-level model, i.e., the owners of the power grid and the natural gas grid. Therefore, the lower-upper model for each path can be expressed as 
\begin{eqnarray}
\label{Eq.OptiTree1}
\begin{aligned}
& {\rm Obj.} \quad (\ref{Eq.LowerObj}) \\ 
& {\rm s.t.} \quad\;\; (\ref{Eq.LowerCon1})-(\ref{Eq.LowerCon33}) \\ 
& \quad\quad\quad{\rm{Given}} \quad v_{l,t}, v_{u,t}, v_{p,t} \quad  \forall l, p, t 
\end{aligned}
\end{eqnarray}
where (\ref{Eq.OptiTree1}) is a mixed integer linear program, which is solved with the CPLEX solver in this research. After solving (\ref{Eq.OptiTree1}), the payoffs of the owners of the natural gas grid and the power grid for each path can be obtained by using (\ref{Eq.UpperEleObj}) and (\ref{Eq.UpperGasObj}). In this case, a finite game tree with payoffs can be established.
\subsection{Game Tree Solution}
Backward induction is a deduction process operating backwards from the end of the game tree to determine a sequence of optimal activities for different players at decision epochs. The optimal decision at each decision node can be expressed as follows.
\begin{align}
\label{Eq.Tree1}
& d_{{h_{i,e}}}^ *  = \arg \mathop {\max }\limits_{{d_{{h_{i,e}}}} \in {D_{{h_{i,e}}}}} {F_i}\left( {\nu ({h_{i,e}}),{d_{{h_{i,e}}}},{\zeta ^ * }(\aleph ({d_{{h_{i,e}}}}))} \right)\quad \forall {h_{i,e}} \\ 
\label{Eq.Tree2}
& {\zeta ^ * }({h_{i,e}}) = \left\{ {d_{{h_{i,e}}}^ * ,{\zeta ^ * }(\aleph (d_{{h_{i,e}}}^ * ))} \right\}\quad \forall {h_{i,e}}  
\end{align}
where (\ref{Eq.Tree1}) shows the optimal decision at the decision node $h_{i,e}$, and $(\nu ({h_{i,e}}),{d_{{h_{i,e}}}},{\zeta ^ * }(\aleph ({d_{{h_{i,e}}}})))$ constructs a path from the initial decision node to a terminal decision node. (\ref{Eq.Tree2}) shows the optimal path from the decision node $h_{i,e}$ to the terminal decision node. By solving the established game tree in sections III.A and III.B with (\ref{Eq.Tree1}) and (\ref{Eq.Tree2}), maintenance plans corresponding to SPNE can be obtained.

\section{Case Studies}
In this section, two test systems are used to verify the proposed model. The first system is a 6-bus power grid with a 4-node natural gas grid, and the second system is a modified IEEE 118-bus power grid with a 20-node natural gas grid. The cases are tested in MATLAB 2017a using {\ct {the CPLEX 12.6 solver}} on computers with 3.1 GHz i5 processors and 8 GB RAMs.
\subsection{6-bus Power Grid with 4-Node Natural Gas Grid}
\subsubsection{Data description}
\begin{figure}[!h]
	\centering
	\includegraphics[width=8cm]{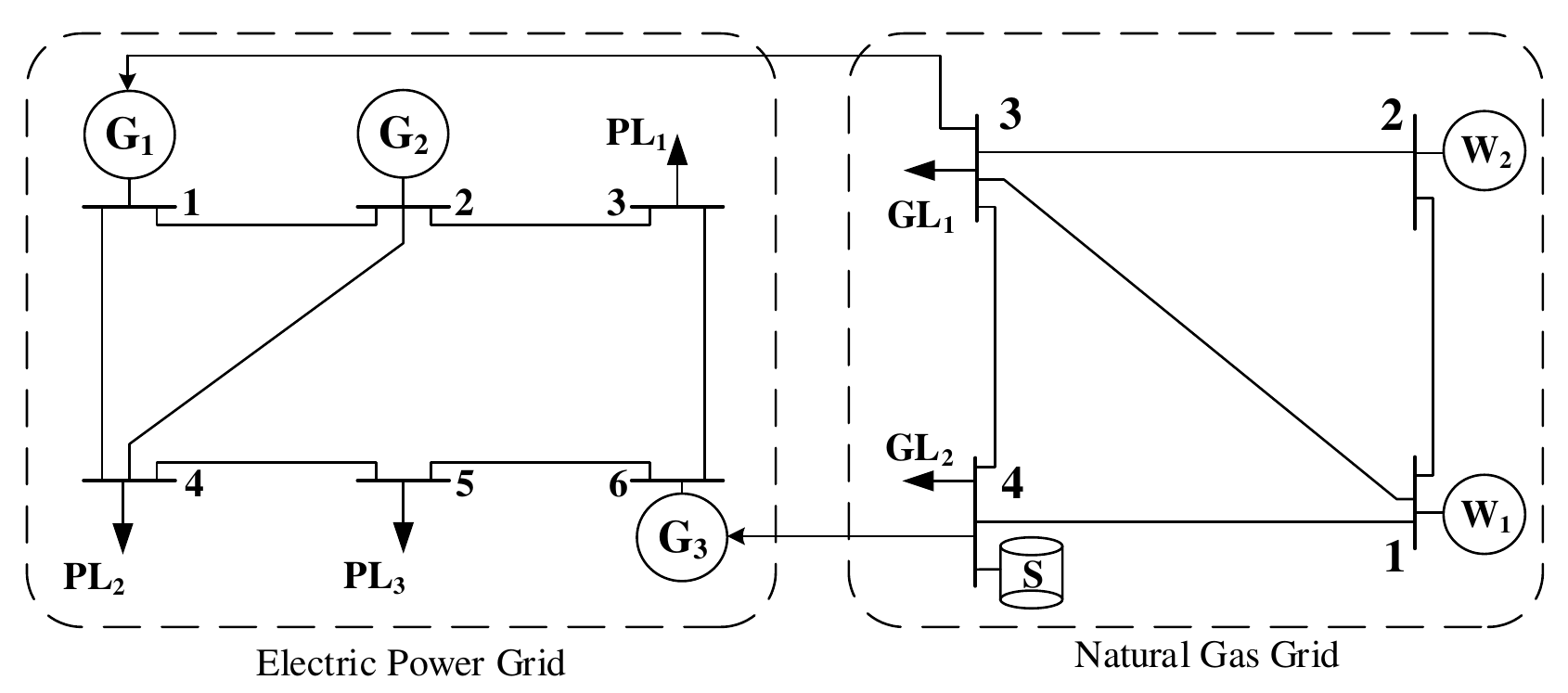}
	\caption{Topology of a 6-bus electric grid with a 4-node gas grid.}
	\label{sixbustest_fig}
\end{figure}

\begin{figure}[!h]
	\centering
	\includegraphics[width=8cm]{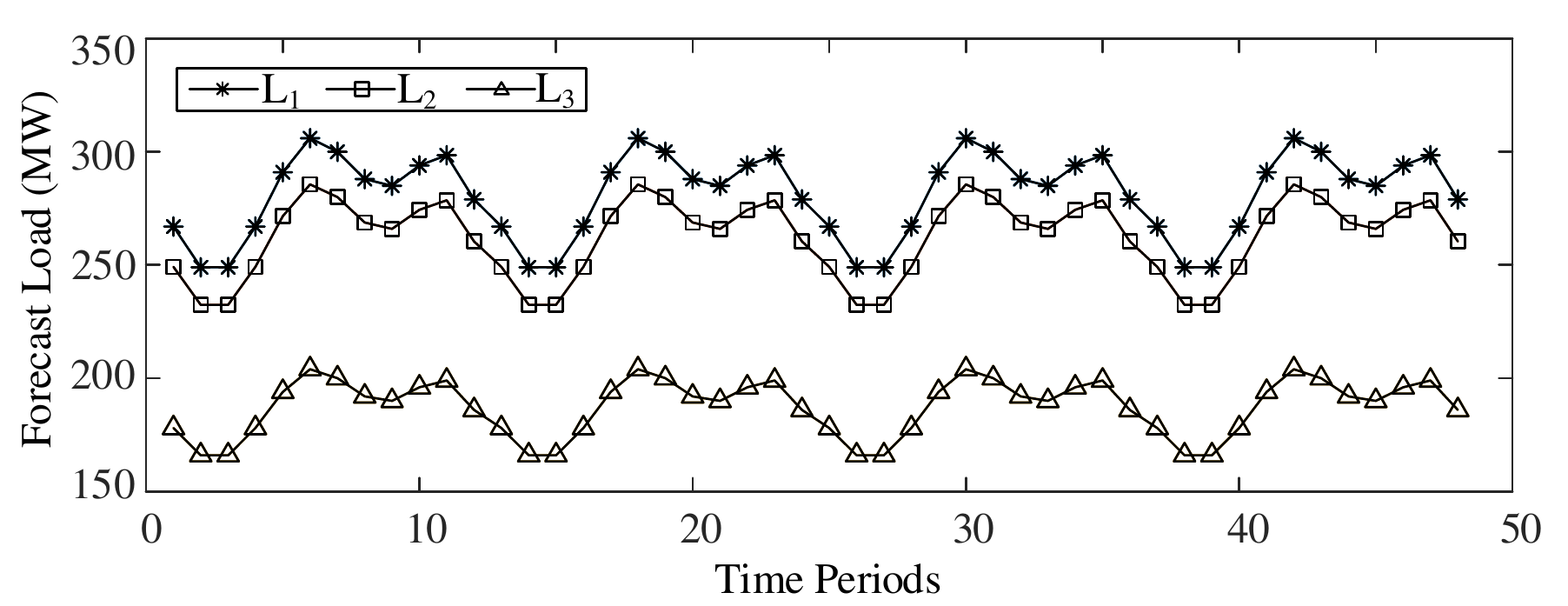}
	\caption{Forecast power load curve.}
	\label{loadcurve_fig}
\end{figure}

The 6-bus power grid and the 4-node natural gas grid refer to \cite{Data1} and \cite{LinerP2}. The integrated system is shown in Fig. \ref{sixbustest_fig}. The forecast power load curve during the given time window is shown in Fig. \ref{loadcurve_fig}.



During the given time window, the lines 1-4 and 3-6 in the power grid will be under maintenance, and maintenance activities need 22 and 24 time intervals, respectively. In each time interval, only one line can be under maintenance. Their maintenance costs are $\$500$ and $\$600$ per time interval. In addition, maintenance on the pipeline 2-3 in the natural gas grid should be implemented with $\$500$ per time interval. The maintenance activities need 31 time periods.

\begin{table}[!h]
	\centering
	\caption{Parameters for Generators}
	\renewcommand{\arraystretch}{1.3}
	\label{table_gene}
	\begin{tabular}{p{3.2cm}<{\centering}p{1.2cm}<{\centering}p{1.2cm}<{\centering}p{1.2cm}<{\centering}}
		\hline
		& $G_1$ & $G_2$ & $G_3$ \\ \hline
		Lower Limits (MW)    & 100   & 80    & 150   \\ 
		Upper Limits (MW)    & 300   & 200   & 350   \\ 
		Ramping Rates (MW/h) & 25    & 20    & 7.5   \\ 
		Fixed Cost (\$)      & 1000  & 1000  & 1000  \\ 
		Linear Cost (\$/MW)  & 20000 & 21000 & 21500 \\ 
		Restart Cost (\$)    & 2000  & 2100  & 2200  \\ 
		Minimum Up Periods   & 4     & 3     & 2     \\ 
		Minimum Down Periods & 2     & 3     & 3     \\ \hline
	\end{tabular}
\end{table}

\begin{table}[!h]
	\centering
	\caption{Possible Maintenance on Line 1-4 ($L_{14}$) and Line 3-6 ($L_{36}$)}
	\renewcommand{\arraystretch}{1.8}
	\label{table.powermain}
	\begin{tabular}{p{0.5cm}<{\centering}p{3cm}<{\centering}|p{0.5cm}<{\centering}p{3cm}<{\centering}}
		\hline
	 No.	& Maintenance Scheduling & No. & Maintenance Scheduling \\ \hline
	$m_{l,1}$  & $L_{14}$:1-22; $L_{36}$:23-46    & $m_{l,7}$     & $L_{14}$:25-46; $L_{36}$:1-24  \\ 
	$m_{l,2}$  & $L_{14}$:1-22; $L_{36}$:24-47    & $m_{l,8}$     & $L_{14}$:26-47; $L_{36}$:1-24   \\ 
	$m_{l,3}$  & $L_{14}$:1-22; $L_{36}$:25-48    & $m_{l,9}$     & $L_{14}$:27-48; $L_{36}$:1-24   \\ 
	$m_{l,4}$  & $L_{14}$:2-23; $L_{36}$:24-47    & $m_{l,10}$    & $L_{14}$:26-47; $L_{36}$:2-25  \\ 
	$m_{l,5}$  & $L_{14}$:2-23; $L_{36}$:25-48    & $m_{l,11}$    & $L_{14}$:27-48; $L_{36}$:2-25 \\ 
	$m_{l,6}$  & $L_{14}$:3-24; $L_{36}$:25-48    & $m_{l,12}$    & $L_{14}$:27-48; $L_{36}$:3-26  \\ \hline
	\end{tabular}
\end{table}

\begin{table}[!h]
	\centering
	\caption{Possible Maintenance on Pipeline 2-3 ($PL_{23}$)}
	\renewcommand{\arraystretch}{1.8}
	\label{table.gasmain}
	\begin{tabular}{p{0.5cm}<{\centering}p{3cm}<{\centering}|p{0.5cm}<{\centering}p{3cm}<{\centering}}
		\hline
		No.	& Maintenance Scheduling & No. & Maintenance Scheduling \\ \hline
		$m_{p,1}$  & $PL_{23}$:1-31    & $m_{p,10}$     & $PL_{23}$:10-40  \\ 
		$m_{p,2}$  & $PL_{23}$:2-32    & $m_{p,11}$     & $PL_{23}$:11-41   \\ 
		$m_{p,3}$  & $PL_{23}$:3-33    & $m_{p,12}$     & $PL_{23}$:12-42  \\ 
		$m_{p,4}$  & $PL_{23}$:4-34    & $m_{p,13}$     & $PL_{23}$:13-43  \\ 
		$m_{p,5}$  & $PL_{23}$:5-35    & $m_{p,14}$     & $PL_{23}$:14-44  \\ 
		$m_{p,6}$  & $PL_{23}$:6-36    & $m_{p,15}$     & $PL_{23}$:15-45  \\ 
		$m_{p,7}$  & $PL_{23}$:7-37    & $m_{p,16}$     & $PL_{23}$:16-46  \\ 
		$m_{p,8}$  & $PL_{23}$:8-38    & $m_{p,17}$     & $PL_{23}$:17-47  \\ 
		$m_{p,9}$  & $PL_{23}$:9-39    & $m_{p,18}$     & $PL_{23}$:18-48  \\ \hline
	\end{tabular}
\end{table}

\subsubsection{Maintenance corresponding to SPNE}
The possible maintenance activities on lines in the power grid are listed in Table \ref{table.powermain}, and the possible maintenance activities on the pipelines in the natural gas grid are listed in Table \ref{table.gasmain}. 

\begin{figure}[!h]
	\centering
	\includegraphics[width=9cm]{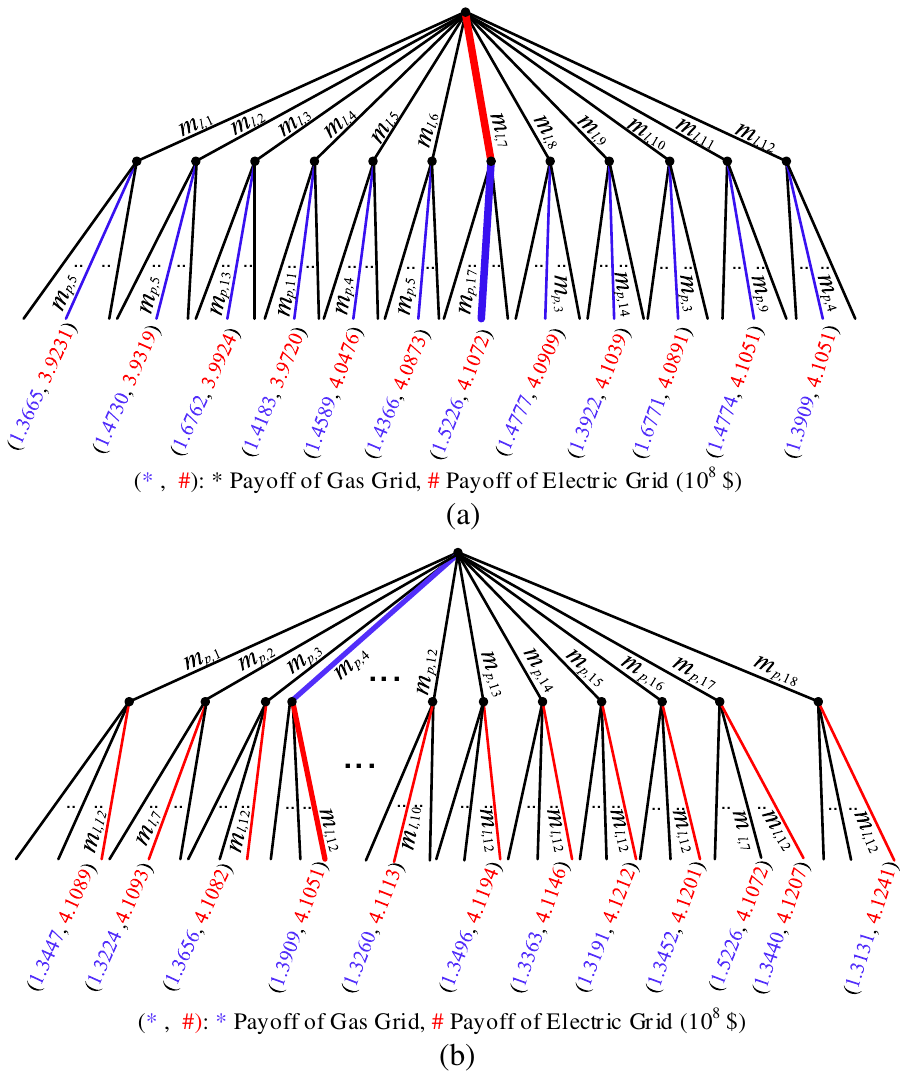}
	\caption{(a) Game tree with payoffs when the owner of the power grid first makes maintenance scheduling. (b) Game tree with payoffs when the owner of the natural gas grid first makes maintenance scheduling.}
	\label{sixbusgametree_fig}
\end{figure}

Fig. \ref{sixbusgametree_fig} (a) shows the game tree when the owner of the power grid first schedules maintenance. A blue path denotes the best response of the owner of the natural gas grid to a maintenance plan scheduled by the owner of the power grid. For example, if the owner of the power grid has the maintenance plan $m_{l,1}$, the best response of the owner of the natural gas to $m_{l,1}$ is $m_{p,5}$. With $m_{l,1}$ and $m_{p,5}$, the payoffs of the owners of the natural gas grid and the power grid are $\$1.3665\times10^8$ and $\$3.9231\times10^8$, respectively. If the owner of the power grid has the maintenance plan $m_{l,3}$, the best response of the owner of the natural gas grid to $m_{l,3}$ is $m_{p,13}$. With $m_{l,3}$ and $m_{p,13}$, their payoffs are $\$1.6762\times10^8$ and $\$3.9924\times10^8$, respectively. Results show that the payoffs are different with diverse maintenance and the responses. Since the two owners expect to maximize their own payoffs, they should take into full account each other's maintenance. The backward induction method is used to obtain the optimal strategies corresponding to SPNE for the two sides. For the scenario in Fig. \ref{sixbusgametree_fig} (a), the optimal strategies corresponding to SPNE is $m_{l,7}$ and $m_{p,17}$. The payoffs for the owners of the power grid and the natural gas grid are $\$4.1072\times10^8$ and $\$1.5226\times10^8$, respectively.

Fig. \ref{sixbusgametree_fig} (b) shows the game tree when the owner of the natural gas grid first schedules maintenance. The red path represents the best response of the owner of the power grid to a maintenance plan scheduled by the owner of the natural gas grid. For example, if the owner of the natural gas grid has the maintenance scheduling $m_{p,1}$, the best response of the owner of the power grid to $m_{p,1}$ is $m_{l,12}$. With the backward induction method, the optimal strategies, which correspond to SPNE, for the two sides are $m_{p,4}$ and $m_{l,12}$. The payoffs of the owners of the power grid and the natural gas grid are $\$4.1051\times10^8$ and $\$1.3909\times10^8$, respectively.

The payoffs corresponding to SPNE in Fig. \ref{sixbusgametree_fig} (b) are both smaller than those in Fig. \ref{sixbusgametree_fig} (a). 
This is because the strategy with a larger payoff of one owner may be excluded from the optimal strategy when the other owner makes the maintenance plan. For example, $m_{l,7}$ and $m_{p,17}$ in Fig. \ref{sixbusgametree_fig} (a) correspond to SPNE when the owner of the power grid first schedules maintenance. In Fig. \ref{sixbusgametree_fig} (b), $m_{l,12}$ rather than $m_{l,7}$ is the best response to $m_{p,17}$ since $\$4.1207\times10^8$ is greater than $\$4.1072\times10^8$.  
Based on the above analysis, we can conclude that different sequences have great influences on the optimal strategies, i.e., SPNE, and the payoffs.

\begin{figure}[!h]
	\centering
	\includegraphics[width=7.5cm]{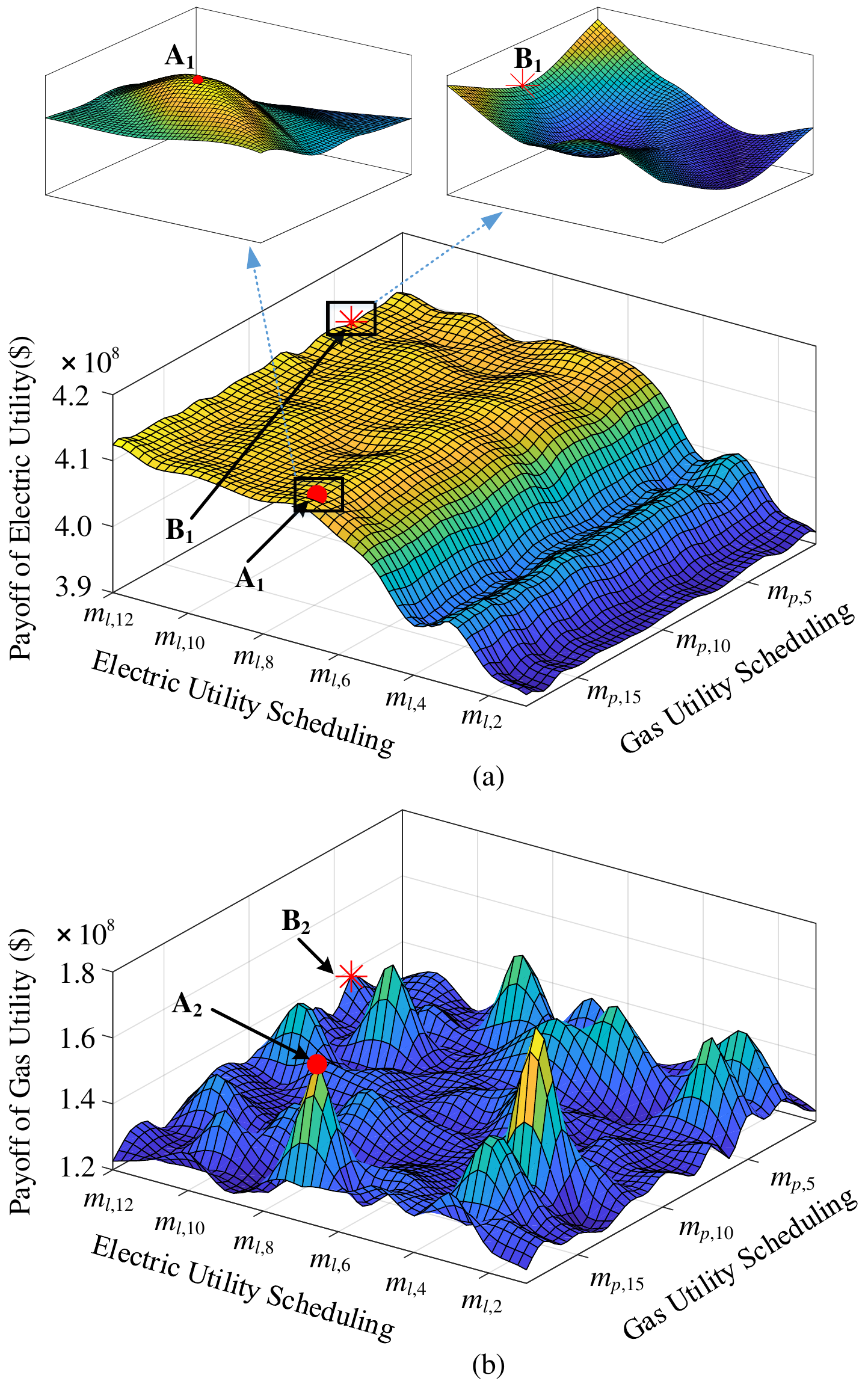}
	\caption{(a) Payoffs of the owner of the power grid under different maintenance plans. (b) Payoffs of the owner of the natural gas grid under different maintenance plans.}
	\label{sixbuspayoff_fig}
\end{figure}

Fig. \ref{sixbuspayoff_fig} (a) and (b) show the fitting surfaces of the payoffs. $A_1$ and $A_2$ represent the payoffs corresponding to SPNE for the power grid and the natural gas grid when the owner of the power grid first schedules maintenance.
$B_1$ and $B_2$ represent the payoffs corresponding to SPNE for the power grid and the natural gas grid when the owner of the natural gas grid first schedules maintenance.
The results show that SPNE may not have the largest payoffs. This is because one strategy with the largest payoff for one owner may be excluded when the other owner makes maintenance plans. Due to this reason, SPNEs are located at different positions of the fitting surfaces of the payoffs, e.g., at the local peaks $A_1$/$A_2$/$B_2$ and the non-local peak $B_1$. To further analyze the position of SPNEs, 400 scenarios with different cost data are implemented, and we get 800 SPNEs in consideration of different sequences of the owners. For these SPNEs, there are 591 SPNEs locating at the local peaks, i.e., about $74\%$ of SPNEs locate at the local peaks.

\begin{figure}[!h]
	\centering
	\includegraphics[width=7.5cm]{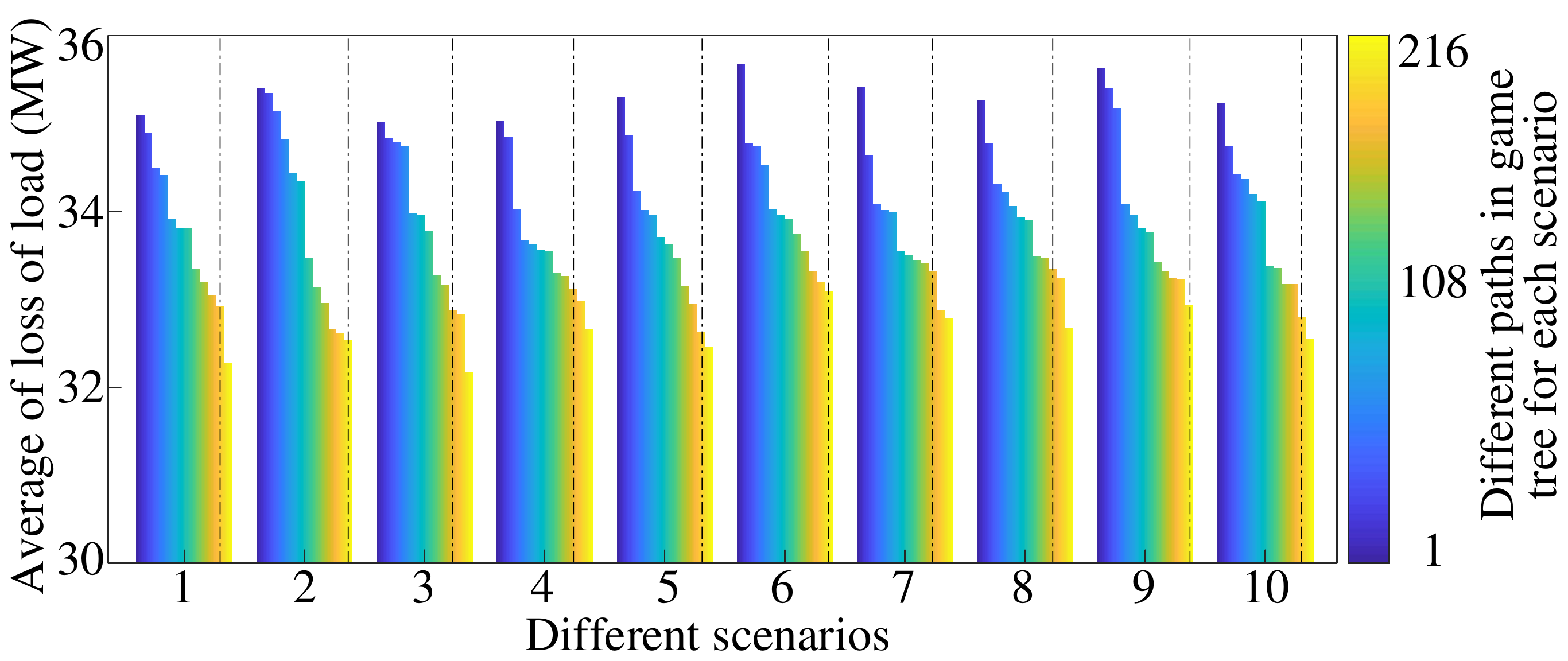}
	\caption{Load shedding of each path in the game tree for different scenarios.}
	\label{LoadShedding_fig}
\end{figure}

Fig. \ref{LoadShedding_fig} shows the average load shedding of each path in the game tree for different scenarios. The x-axis represents the scenario, and the y-axis represents the average load shedding. Different color bars in each scenario represent the different paths in the game tree. Each black dash-dotted line marks the path corresponding to SPNE in each scenario. Some scenarios (e.g., the scenarios 2, 6 and 9) have SPNEs with the minimum loss of load, and the others (e.g., the scenarios 1, 3, 4, 5, 7, 8, and 10) have SPNEs with the non-minimum loss of load. One critical reason for this is that the strategy with the minimum loss of load may be excluded when the other player makes maintenance plans.

\subsubsection{Influences of maintenance duration}  
\begin{figure}[!h]
	\centering
	\includegraphics[width=7cm]{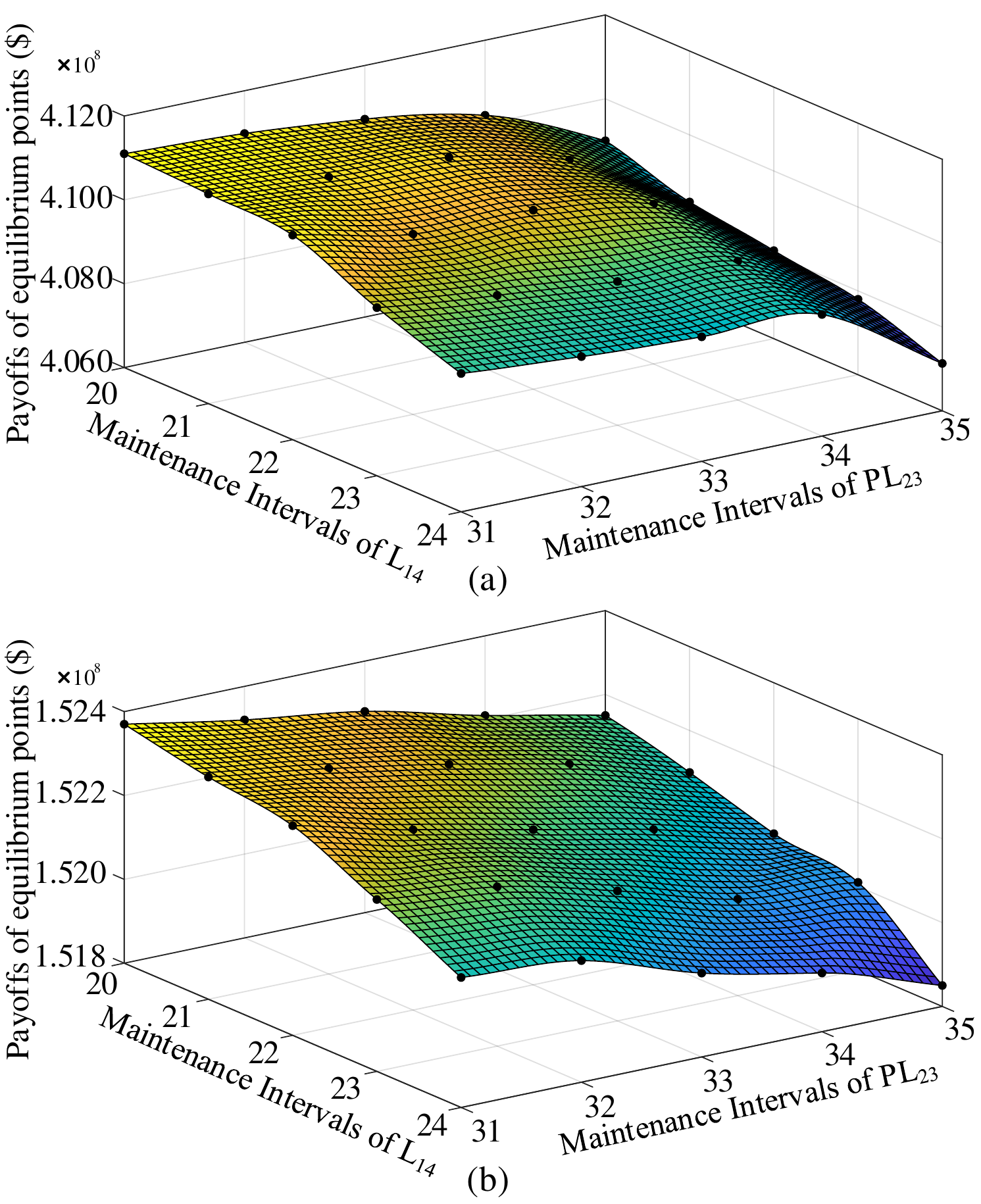}
	\caption{Payoffs of the equilibrium points with different maintenance durations for the owner of the power grid (a) and the owner of the natural gas grid (b), respectively.}
	\label{DifferentDuration_fig}
\end{figure}

Fig. \ref{DifferentDuration_fig} (a) and (b) show the payoffs corresponding to  SPNEs with different maintenance durations. The x-axis and the y-axis represent the maintenance durations of the line $L_{14}$ and the pipeline $PL_{23}$, respectively. The z-axis represents the payoffs. The results show that the payoffs corresponding to SPNEs gradually increase with the shorter maintenance durations. We can interpret this pattern qualitatively as follows. 
\begin{itemize}
	\item When the maintenance duration of the line $L_{14}$ decreases, the maintenance cost for the line can be reduced. Furthermore, the loss of load could be reduced due to a larger power transport capability, and in consequence the generation revenue increases.
	\item  When the maintenance duration of the pipeline $PL_{23}$ decreases, the maintenance cost for the pipeline can be reduced, and the loss of gas load could be reduced. In addition, the generation from the gas-fired units can increase due to the larger gas transport capacity, and in consequence the loss of load can be reduced.  
\end{itemize}

Based on the results and the analysis, it is concluded that the owners of the natural gas grid and the power grid should shorten the maintenance durations when they would like to obtain SPNE with higher payoffs.

\subsubsection{Influences of different piecewise linear functions}
\begin{figure}[!h]
	\centering
	\includegraphics[width=7.5cm]{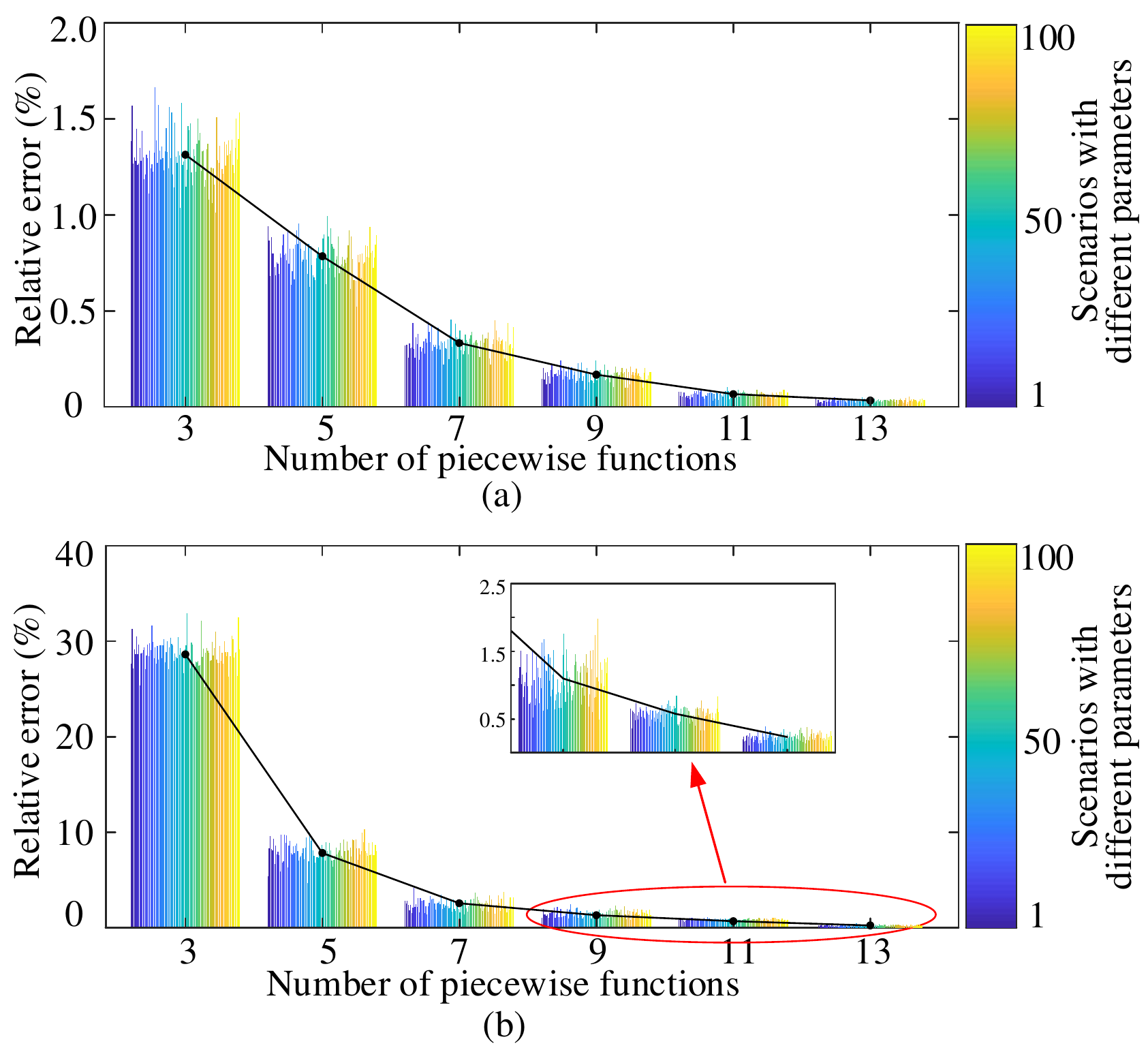}
	\caption{Relative errors of the payoffs corresponding to the equilibrium points for the owner of the power grid (a) and the owner of the natural gas grid (b), respectively.}
	\label{Error_piecewise_fig}
\end{figure}

Fig. \ref{Error_piecewise_fig} shows the relative errors of the payoffs corresponding to SPNE for the power grid (a) and the natural gas grid (b) when using different piecewise functions. The x-axis and the y-axis represent the number of piecewise functions and the relative error, respectively. The different color bars represent scenarios with different cost data. The results with fifteen piecewise functions are used as reference values. The results are more accurate with more piecewise functions, and nine piecewise functions can achieve the results with high accuracy.

Furthermore, the payoffs of the owner of the natural gas grid are more sensitive to the number of piecewise functions compared to the payoffs of the owner of the power grid. For example, the average relative error of the payoff is about 1.4\% for the power grid, and about 28\% for the natural gas grid when using three piecewise functions. When using nine piecewise functions, the average relative error of the payoff is about 0.25\% for the power grid, and about 1\% for the natural gas grid.     

\subsubsection{Influences of maintenance on different components}
Maintenance on generators may need to be scheduled associated with maintenance on lines during the given time window. To compare the results, the line $L_{36}$ in the above case is replaced by the generator $G_{2}$ that needs maintenance, and the other settings remain unchanged. We compare the results when the owner of the power grid first makes maintenance scheduling. Table \ref{table_gasunits_payoffs} shows the comparison results. The scenario $S_1$ is the original case, and the scenario $S_2$ is the updated case with the generator $G_{2}$ to be under maintenance.

\begin{table}[!h]
	\centering
	\renewcommand{\arraystretch}{1.5}
	\caption{Payoffs Corresponding to Equilibrium Points}
	\label{table_gasunits_payoffs}
	\begin{tabular}{p{1.8cm}<{\centering}p{2.8cm}<{\centering}p{2.8cm}<{\centering}}
		\hline
		\multirow{2}{*}{Scenario} & \multicolumn{2}{c}{Payoffs ($\$10^8$)} \\ \cline{2-3} 
		& Power Grid                                & Natural Gas Grid                               \\ \hline
		$S_1$        & 4.1072                                          & 1.5226                                    \\ \hline
		$S_2$        & 3.4550                                          & 1.7761                                    \\ \hline
	\end{tabular}
\end{table}

The payoff of the scenario $S_2$ corresponding to the equilibrium point for the power grid ($\$3.4550 \times10^8$) is smaller than that of the scenario $S_1$ ($\$4.1072 \times10^8$). The main reason for this is that the loss of load increases when the generator $G_{2}$ is under maintenance.

The payoff of the scenario $S_2$ corresponding to the equilibrium point for the natural gas grid ($\$1.7761 \times10^8$) is larger than that of the scenario $S_1$ ($\$1.5226 \times10^8$). This is because the gas-fired units $G_1$ and $G_3$ need more natural gas for electric power to reduce the loss of load when the generator $G_{2}$ is under maintenance, and the increased natural gas production increases the revenues of the natural gas grid.

\subsection{118-bus electric grid with 20-node gas grid}
\subsubsection{Data description}
The data for the power grid and the natural gas grid refer to \cite{Data2} and \cite{Data3}, respectively. Table \ref{table_gasunit} shows the gas-fired units and the corresponding gas nodes in the natural gas grid. During the given time window (48 time intervals), three transmission lines 69-47, 69-70 and 69-75 will be maintained with 24 time intervals, and two pipelines 10-14 and 10-11 will be maintained with 23 time intervals. The maximum numbers of lines and pipelines that can be maintained in one period are 2 and 1, respectively. 

\begin{table}[!h]
	\centering
	\renewcommand{\arraystretch}{1.4}
	\caption{Number of Gas-fired Units}
	\label{table_gasunit}
	\begin{tabular}{p{2.5cm}<{\centering}p{2.5cm}<{\centering}p{2.5cm}<{\centering}}
		\hline
		Gas-fired Unit No. & Power Bus & Gas Node \\ \hline
		1                  & 10        & 16       \\ 
		2                  & 12        & 15       \\ 
		3                  & 54        & 12       \\ 
		4                  & 59        & 18       \\ 
		5                  & 87        & 5        \\ 
		6                  & 103       & 19       \\ \hline
	\end{tabular}
\end{table}

\subsubsection{Simulation Results}
This section shows the maintenance scheduling (the equilibrium) for the two owners with different decision-making sequences. Table \ref{table.118_1} shows the maintenance scheduling and the payoffs when the owner of the power grid first makes the maintenance plan. Table \ref{table.118_2} shows the maintenance scheduling and the payoffs when the owner of natural gas grid first makes the maintenance plan. The results show that different decision making sequences result in different SPNEs, and therefore cause different maintenance plans. 

\begin{table}[!h]
	\centering
	\caption{Payoffs with Power Grid's Owner First Making Maintenance Plans}
	\renewcommand{\arraystretch}{1.8}
	\label{table.118_1}
	\begin{tabular}{p{1.0cm}<{\centering}p{2.2cm}<{\centering}|p{1.0cm}<{\centering}p{2.4cm}<{\centering}}
		\hline
		\multicolumn{2}{c|}{Power Grid} & \multicolumn{2}{c}{Natural Gas Grid} \\ 
		Payoff ($\$10^8$)  & Maintenance Scheduling  & Payoff ($\$10^8$)   & Maintenance Scheduling   \\ \hline
		$52.1343$            & $L_{69-47}$: $22-45$ $L_{69-70}$: $25-48$ $L_{69-75}$: $1-24$                       & $6.6483$             & $PL_{10-14}$: $2-24$ $PL_{10-11}$: $26-48$                       \\ \hline
	\end{tabular}
\end{table}

\begin{table}[!h]
	\centering
	\caption{Payoffs with Natural Gas Grid's Owner First Making Maintenance Plans}
	\renewcommand{\arraystretch}{1.8}
	\label{table.118_2}
	\begin{tabular}{p{1.0cm}<{\centering}p{2.2cm}<{\centering}|p{1.0cm}<{\centering}p{2.4cm}<{\centering}}
		\hline
		\multicolumn{2}{c|}{Power Grid} & \multicolumn{2}{c}{Natural Gas Grid} \\ 
		Payoff ($\$10^8$)  & Maintenance Scheduling  & Payoff ($\$10^8$)   & Maintenance Scheduling   \\ \hline
		$52.1321$            & $L_{69-47}$: $25-48$ $L_{69-70}$: $1-24$ $L_{69-75}$: $21-41$                       & $6.7062$             & $PL_{10-14}$: $2-24$ $PL_{10-11}$: $25-47$                       \\ \hline
	\end{tabular}
\end{table}

We shorten the maintenance durations of the line $L_{69-47}$ and the pipeline $L_{10-14}$ to $20$ and $19$, respectively. Table \ref{table_118_new} shows the payoffs corresponding to the equilibrium points. When the owner of the power grid first makes the schedule, the payoffs for the owners of the power grid and the natural gas grid are $\$52.4254\times10^8$ and $\$6.8107\times10^8$, respectively. They are larger than the scenario in Table \ref{table.118_1}. When the owner of the natural gas grid first makes the schedule, the payoffs for the owners of the power grid and the natural gas grid are $\$52.4172\times10^8$ and $\$6.8211\times10^8$, respectively. They are larger than the scenario in Table \ref{table.118_2}. The results show that the payoffs corresponding to the equilibrium point increase when the owners shorten the maintenance durations. This conclusion is consistent with that in Section IV.A.3.

\begin{table}[!h]
	\centering
	\renewcommand{\arraystretch}{1.4}
	\caption{Payoffs Corresponding to SPNE under Different Maintenance Durations}
	\label{table_118_new}
	\begin{tabular}{ccc}
		\hline
		\multirow{2}{*}{Sequence of making schedules} & \multicolumn{2}{c}{Payoffs ($\$10^8$)} \\ \cline{2-3} 
		& Power Grid  & Natural Gas Grid \\ \hline
		Power grid first                        & $52.4254$            & $6.8107$      \\ \hline
		Natural gas grid first                            & $52.4172$            & $6.8211$      \\ \hline
	\end{tabular} 
\end{table}

\subsection{Discussion}
In this paper, the natural gas grid and the power grid are coupled via natural gas-fired units. In addition to the gas-fired units, there are some other ways such as energy hubs for receiving, sending, converting and storing different types of energy, and these energy hubs are also critical to the natural gas system and the power system. There are many research studies on the steady-state models of the energy hubs \cite{hub1, hub2, hub3}. When an energy hub is included in the natural gas grid, its steady-state model can be expressed as $\bf{O}=\bf{C} \times \bf{I}$, in which $\bf{I}$ and $\bf{O}$ are vectors of energy at input and output ports, and $\bf{C}$ is the coupling matrix that describes the conversion of energy from the input to the output. The elements in $\bf{C}$ represent the converter efficiencies and the hub internal topology \cite{hub4}. This steady-state model has a similar form as the relation between natural gas and electric power, i.e., the constraint (\ref{Eq.LowerCon28}). Therefore, the the model of the energy hub can be easily integrated in the proposed model.

Currently, we mainly consider maintenance on power lines, gas pipelines and the generating units. We would like to emphasize that the proposed dynamic game-based model is suitable for other components in the system. The models for the other components can be established by analogy with lines, pipelines and generating units. For example, the maintenance model for a transformer can be established by analogy with the maintenance model for the line, and the maintenance model for a gas well can be established by analogy with the maintenance model for the generating unit.

When modeling the natural gas system in this paper, a piecewise linear approximation approach is used to transform the nonconvex and nonlinear Weymouth gas flow equations into a series of piecewise linear functions. Since maintenance scheduling is a planning problem, appropriate piecewise linear functions can satisfy the accuracy requirements. To further improve the accuracy, some state-of-the-art methods based on the second-order conic relaxation \cite{Convex1, Convex5} can be used to deal with the Weymouth equations. When using the second-order conic relaxation, we just need to replace the piecewise linear model for Weymouth equations by the second-order conic relaxation, and the optimization model in each game tree becomes a mixed integer nonlinear programming (MINLP) model.

{\ct{For a large scale system, it may need a long computation time to obtain the maintenance plans. However, maintenance scheduling belongs to a planning topic, it does not require online calculation. In addition, one critical point for power systems and natural gas grids is to ensure a high reliability level. To this end, the number of transmission lines/pipelines that can be out of service in the same time period is limited \cite{EPRI11}. Meanwhile, the number of paths is also constrained by the continuous minimum periods for maintenance activities. Considering the above practical points, i.e., the constraints (2)-(8), the number of paths is limited, particularly when the number of all transmission lines/pipelines to be under maintenance with a given time window is limited. Furthermore, many simulation platforms such as Rescale \cite{Rescale11}, IBM Could \cite{IBMCloud11} and Nimbix \cite{Nimbix11} with high performance computing (HPC), can be used to solve the proposed model efficiently.}}

Some techniques can be used to reduce the computation time. First, parallel computing can be employed when constructing the game tree since each path in the game tree are independent \cite{Parallel1}. Second, some methods, e.g., the parallel branch \& bound method \cite{BBP1} and {\ct{the modified Benders decomposition method \cite{Secondrefer}}}, can be used to accelerate the convergence of the solution for each path. 

\section{Conclusion}
This paper proposes a dynamic game-based maintenance scheduling mechanism for natural gas and power grids by using a bilevel approach. In the upper-level model, the different owners of the power grid and the natural gas grid schedule maintenance with the objective to maximize their own revenues. This is modeled as a dynamic game problem, which is solved by the backward induction algorithm. In the lower-level model, the operating point under the scheduled maintenance is calculated. This is modeled as a mixed integer linear programming problem. For the model of natural gas grids, the piecewise linear approximation associated with the big-M approach is used to transform the original nonlinear model into a mixed integer linear model. 

The proposed model was validated by two test systems. The major findings are as follows. 1) The subgame perfect Nash equilibrium for the power grid and the natural gas grid may not have the maximum payoffs, however, it has a high probability to locate a local peak of the payoff surface. 2) The maintenance duration has a great impact on the payoff corresponding to the subgame perfect Nash equilibrium. When the owners of the power grid and the natural gas grid expect to obtain an equilibrium point with the higher payoffs, it is suggested that they should shorten the maintenance durations. 3) The payoff corresponding to the equilibrium for the owner of the natural gas grid has a high sensitivity to the number of piecewise functions compared to the payoff corresponding to the equilibrium for the owner of the power grid.

This work can be extended to consider the integration of renewables. Furthermore, some state-of-the-art methods based on the second-order conic relaxation can be used to deal with the Weymouth equations, and the optimization model in each game tree can be established as a mixed integer nonlinear programming (MINLP) model.




\bibliographystyle{IEEEtran}
\bibliography{IEEEabrv,RefDatabase}

\end{document}